\documentclass[twoside,a4paper]{article}
\input amssym.def
\input amssym.tex
\font\tencmmib=cmmib10
\skewchar\tencmmib='177
\font\sevencmmib=cmmib7
\skewchar\sevencmmib='177
\font\fivecmmib=cmmib5
\skewchar\fivecmmib='177
\newfam\cmmibfam
\textfont\cmmibfam=\tencmmib
\scriptfont\cmmibfam=\sevencmmib
\scriptscriptfont\cmmibfam=\fivecmmib
\font\tencmbsy=cmbsy10
\skewchar\tencmbsy='60
\font\sevencmbsy=cmbsy7
\skewchar\sevencmbsy='60
\font\fivecmbsy=cmbsy5
\skewchar\fivecmbsy='60
\newfam\cmbsyfam
\textfont\cmbsyfam=\tencmbsy
\scriptfont\cmbsyfam=\sevencmbsy
\scriptscriptfont\cmbsyfam=\fivecmbsy
\newtheorem{Definition}{Definition}
\newtheorem{Proposition}{Proposition}
\newtheorem{Theorem}{Theorem}
\newtheorem{Lemma}{Lemma}
\newtheorem{Corollary}{Corollary}

\newtheorem{Remark}{Remark}
\title{Locally nearly spherical surfaces are almost-positively $c$-curved\footnote{2000 {\it Mathematics Subject Classification}: Primary 53C21; Secondary 49N60; {\it Keywords}: Monge problem, compact surfaces, nearly spherical, positive $c$-curvature, stability}}
\author{Ph. Delano\" e\thanks{Supported by the CNRS}\hspace{.4em} and Yuxin Ge}
\date{}
\begin{document}
\maketitle
\thispagestyle{empty}
\begin{abstract}
The $c$-curvature of a complete surface with Gauss curvature close to 1 in $C^2$ norm is almost-positive (in the sense of Kim--McCann). Our proof goes by a careful
case by case analysis combined with perturbation arguments from the constant curvature case, keeping track of an estimate on the closeness curvature
condition.
\end{abstract}
\section{Introduction and main results}
Monge's problem, in optimal transport theory, goes back to \cite{mon}. In its general formulation, one is looking for an optimal map
$f:(M,\mu)\to(\overline{M},\bar{\mu})$ between two Polish probability spaces. The optimality criterion consists in minimizing the total cost
functional%
\def\ds{\displaystyle}%
$\ds\int_M c(x,f(x))\ d\mu(x)$ among measurable maps which push $\mu$ to $\bar{\mu}$, where the cost function $c:M\times\overline{M}\to {\Bbb
R}\cup\{+\infty\}$ is given lower semi-continuous with some additional properties (see
\textit{e.g.}
\cite{vil} and references therein). In the emblematic case of the Brenier--McCann cost function: $\ds M=\overline{M}, c=\frac{1}{2}d^2,$ where
$M$ stands for a complete Riemannian manifold with associated distance function $d$, this problem was solved under mild assumptions on the given
probability measures $\mu$ and $\bar{\mu}$ \cite{bre,mcc}. In that case, the optimal map must read %
\def\grad{\mathop{\rm grad}\nolimits}%
$f=\exp(\grad u)$ for some $c$-convex potential
function $u$ such that the pushing condition $f_{\#}\mu=\bar{\mu}$ becomes a partial differential equation of Monge--Amp\`ere type satisfied by $u$ in a weak sense.
Neil Trudinger and his co-workers observed that a similar solution scheme exists for a class of more general cost functions $c$ for which,
given smooth data, they analyzed the smoothness of the corresponding potential function $u$ \cite{mtw}. For the purpose of a one-sided interior
estimate on an expression of second order (in $u$), they were lead to formulate a fourth-order two-points condition on the cost function $c$,
called (A3S) condition. A weak form of the latter, called (A3W), was proved necessary (for the smoothness of $u$) by Loeper \cite{loe};
in particular, in the Brenier--McCann case, he interpreted (A3W) read on the diagonal of $M\times M$ as the non-negativity of the sectional
curvature of $M$. Lately, still with $c=\frac{1}{2}d^2$, C\'edric Villani and his co-workers were able to relate some variants of (A3S), checked stable at round spheres under $C^4$ small deformations of the standard round metric, with the
convexity of the tangential domain of injectivity of the exponential map \cite{lovi,firi,firivi2}. However, the very
geometrical status of the fourth-order expression (in $c$) occuring in condition (A3S) was not understood untill Kim and McCann interpreted it
\cite{kimac1} as a genuine, though quite special, curvature expression arising on the product manifold $M\times\overline{M}$ endowed with the
pseudo-Riemannian  metric: $\ds h=-\frac{1}{2}\frac{\partial^2c}{\partial x^i\partial \bar{x}^j}(dx^i\otimes d\bar{x}^j + d\bar{x}^j \otimes
dx^i)$. They also defined an extended version of (A3S), stronger than (A3W), called non-negative cross-curvature condition (NNCC, for short\footnote{also used below as an abbreviation for non-negatively cross-curved}) and proved that it is
stable under Cartesian product\footnote{unlike (A3S), or even (A3W) alone as soon as a factor is not NNCC} as well as, in the Brenier--McCann case, under Riemannian submersion \cite{kimac2}. Actually, in that case, they defined a stronger condition called almost-positive cross-curvature condition (APCC, for short\footnote{also used below as an abbreviation for almost-positively cross-curved}) also shown stable under Riemannian submersion \cite{kimac2}. So, with $c=\frac{1}{2}d^2$, the stability of NNCC (resp. APCC) under products (resp. submersions) enables to
construct new NNCC (resp. APCC) examples out of known ones -- like the standard sphere \cite{kimac2}. In the present paper, we will prove the stability
of APCC at the standard 2-sphere; specifically, we will check the APCC condition for $c=\frac{1}{2}d^2$ on a complete surface with Gauss curvature $C^2$ close to a positive constant. This result complements the stability one of \cite{firi} on the 2-sphere as well as an unstated one (stability of APCC at the standard $n$-sphere, near conjugacy, see Remark \ref{rem-firivi2} below) obtained in the course of a proof in \cite{firivi2}. Here, let us point out that our paper is drawn from an initial 44 pages draft sent by Ge to Figalli in January 2009, thus independent from the papers \cite{firivi1,firivi2} first circulated in July 2009; in particular, our analysis of the APCC property near conjugacy (Section \ref{sec-conj} below) departs from that of \cite{firivi2}.\\
\def\cut{\mathop{\rm Cut}\nolimits}%
\def\sec{\mathop{\rm Sect}\nolimits}%

In order to state our result, let us first recall some definitions, restricting to connected complete Riemannian manifolds $M=\overline{M}$ with
the cost function
$c=\frac{1}{2}d^2$ defined on $M\times\overline{M}\ \setminus\cut$, where Cut stands for the cut locus. Using the aforementioned
pseudo-Riemannian metric
$h$ on
$M\times\overline{M}$ and setting $\sec_h$ for its sectional curvature tensor viewed as a field of quadratic forms on
$\bigwedge^2T(M\times\overline{M})$, for each
$(m,\overline{m})\in M\times\overline{M}\ \setminus\cut$ and each $(\xi,\bar{\xi})\in T_mM\times T_{\overline{m}}\overline{M}$, the
associated
\textit{cross-curvature} is defined by \cite{kimac1}:%
\def\cro{\mathop{\rm cross}\nolimits}%
$$\cro_{(m,\overline{m})}(\xi,\bar{\xi}):=\sec_h[(\xi\oplus 0)\wedge(0\oplus \bar{\xi})].$$
Kim and McCann observed that it must vanish for some choice of $(\xi,\bar{\xi})$ \cite{kimac2}. If it is identically non-negative, we say that
the manifold
$M$ is NNCC. For instance, the standard
$n$-sphere is NNCC
\cite{kimac2,firi} and if a manifold $M$ endowed with a Riemannian metric $g$ is so, its sectional curvature tensor $\sec_g$ must be
non-negative because, at
$\overline{m}=m$, we have:
$\ds\cro_{(m,m)}(\xi,\bar{\xi})\equiv \frac{4}{3}\sec_g(\xi\wedge\bar{\xi})$, as first observed by
Loeper \cite{loe}.
Pulling back by the exponential map:%
\def\nocut{\mathop{\rm NoCut}\nolimits}%
\def\noconj{\mathop{\rm NoConj}\nolimits}%
$$(m,V)\in\nocut\longrightarrow (m,\exp_m(V))\in M\times\overline{M}³\setminus\cut\ ,$$
where $\nocut$ is the domain of $TM$ defined by:
$$\nocut:=\left\{(m,V)\in TM,\ \forall t\in[0,1], \exp_m(tV)\notin\cut_m \right\}$$
(and $\cut_m$, the cut locus of the point $m$), Trudinger \textit{et al} noted \cite[p.164]{mtw}
that one identically recovers $\frac{1}{2}\cro_{(m,\overline{m})}(\xi,\bar{\xi})$ at $\overline{m}=\exp_m(V)$ with
$(m,V)\in\nocut$ and
$\bar{\xi}=d(\exp_m)(V)(\nu)$, by calculating the quantity:
\begin{equation}\label{c-curvdef-eq}
{\cal C}(m,V)(\xi,\nu):=-\frac{D^2}{d\lambda^2}\left[A(m,V+\lambda\nu)(\xi) \right]_{\lambda=0}
\end{equation}
where: $A(m,V)(\xi):=\nabla d[p\to c(p,\exp_m(V))]_{p=m}(\xi,\xi)$ with $\nabla$ the Levi--Civita connection of the Riemannian metric $g$ and
where $D$ stands for the canonical flat connection of $T_mM$. In \cite{delge}, we
performed a stepwise calculation of $A(m,V)(\xi)$ and its first and second derivatives with respect to $V$, in a Fermi chart along the geodesic
$t\in[0,1]\to\exp_m(tV)\in M$. This calculation just requires that $(m,V)$ belong to $\noconj$, denoting so the domain\footnote{as well-known
\cite{carm,cheb}, $\noconj$ is the maximum rank domain for the exponential map which contains $\nocut$} of $TM$ which consists of tangent vectors $(m,W)\in TM$ such that the geodesic segment $t\in[0,1]\to\exp_m(tW)$
contains no conjugate points, a
fact conceptualized in \cite{firi} using the Hamiltonian flow (see also \cite{firivi2}).
Neil Trudinger suggested that one calls the
quantity
$\ds{\cal C}(m,V)(\xi,\nu)$ defined by (\ref{c-curvdef-eq}), now with $(m,V)\in\noconj$, the
$c$-curvature\footnote{somewhat consistently with the $c$-segment denomination used in \cite{mtw}; we will use this short denomination, instead of 'extended MTW tensor' as in \cite{lovi,firivi1,firivi2} or 'Ma--Trudinger--Wang curvature' as in \cite{firivi1,leemcc}, since further names could be associated to the birth of this conceptual object, anyhow} of $M$ at $(m,V,\xi,\nu)$. It is known to vanish if ${\rm rank}(V,\xi,\nu)\leq 1$ \cite{delge,kimac2}. Now, the definition given in \cite{kimac2} of an APCC (resp. NNCC) manifold reads in terms of the $c$-curvature as follows:
\begin{Definition}\label{apcc-def}
Let $M$ be a connected complete Riemannian manifold with cost function $c=\frac{1}{2}d^2$. We say that $M$ is non-negatively $c$-curved, or NNCC, if $\ds{\cal C}(m,V)(\xi,\nu)\geq 0$ for each $(m,V)\in\noconj$ and each couple $(\xi,\nu)$ in $T_mM$. If $M$ is NNCC and such that: $\ds{\cal C}(m,V)(\xi,\nu)=0$ if and only if the span of $(V,\xi,\nu)$ has dimension at most 1, we call it almost-positively $c$-curved, or APCC.
\end{Definition}
\begin{Remark}\label{rem-firivi2}
{\rm An intermediate (unquantified) result of \cite{firivi2}, unstated there as such, obtained \textit{via} the square completion of a huge expression, goes as follows:}\\
If $M$ is the $n$-sphere endowed with a Riemannian metric $C^4$ close to the standard one and if $(m,V)\in\noconj$ lies close enough to the boundary of $\noconj$, then $(\xi,\nu)\mapsto {\cal C}(m,V)(\xi,\nu)$ satisfies the APCC property on $T_mM\times T_mM$.
\end{Remark}
Let us call, for short, a \textit{surface} any smooth compact connected 2-dimensional Riemannian manifold without boundary. We aim at the
following result:
\begin{Theorem}\label{th-main}
Let $S$ be a surface with Gauss curvature $K$ such that $\min_S K=1$. There exists a small universal constant $\eta >0$ such that, if
$\ds |K-1|_{C^2(S)}\leq\eta$, then $S$ is APCC.
\end{Theorem}
Here, the $C^2$ norm of a function $f:S\to{\Bbb R}$ is defined (using the Riemannian norm $|.|$ on tensors) by:
$\ds |f|_{C^2(S)}:=\sup_S|f|+\sup_S |df|+\sup_S |\nabla df|$.\\
The result is proved in \cite{kimac2} with $\eta=0$ (constant curvature case, see also \cite{firi} for NNCC) and in \cite{firivi2} with $(m,V)\in\noconj$ lying close enough to the boundary of $\noconj$ (with no quantified estimates, though). If $V=0$, the result is obvious
(due to the cross-curvature interpretation when
$m=\overline{m}$), so we will assume
$V\not=0$ with no loss of generality.
\begin{Remark}\label{comparison-rk}
{\rm Let
$${\rm D}_c=\sup \left\{ |V|_m,\ (m,V)\in\noconj\right\}$$
be the diameter of conjugacy of $S$. Since $K\geq 1$, the Bonnet--Myers theorem \cite{bon,carm,cheb,mil} implies: ${\rm D}_c\leq \pi\,$; in
particular, the diameter of
$S$ must be at most equal to $\pi$.}
\end{Remark}

Actually, we will prove a stronger result, namely:
\begin{Theorem}\label{strong-c-posit-th}
Let $S$ be a surface with $\min_S K=1$. There exists small universal positive constants $\eta,\varsigma$ such that, if $\ds |K-1|_{C^2(S)}\leq\eta$, for any $(m,V)\in \noconj$ and any couple $(\xi,\nu)$ of unit vectors in $T_mS$, the following inequality holds:
\begin{equation}\label{strong-c-posit-ineq}
{\cal C}(m,V)(\xi,\nu) \geq \varsigma\ {\cal A}_2(m,V,\xi,\nu)\ ,
\end{equation}
where ${\cal A}_2(m,V,\xi,\nu)$ stands for the sum of the squared areas of the parallelograms repectively defined in $T_mS$ by the couples $(\xi,\nu),(V,\xi),(V,\nu),$ in other words:
$${\cal A}_2(m,V,\xi,\nu)=|\xi |^2 |\nu |^2 - g(\xi,\nu )^2+|V |^2 |\xi |^2 - g(V ,\xi )^2+| V|^2 | \nu|^2 - g(V ,\nu )^2.$$
\end{Theorem}

The outline of the paper essentially coincides with that of the proof. We present a quick derivation of the
$c$-curvature expression in Section
\ref{sec-express}
and related perturbative estimates for that expression, based on the assumption that the
$C^2$ norm of
$(K-1)$ is small, in Section \ref{sec-perturb}. Using the latter, we prove successively Theorem \ref{strong-c-posit-th} under the additional assumption that
the point $\exp_m(V)$ lies, either near the first conjugate point $m^*$ of $m$ along the geodesic $t\in{\Bbb R}^+\to \exp_m(tV)\in S$ (Section \ref{sec-conj}), or near
$m$ (Section
\ref{sec-zero}), or in-between (Section
\ref{sec-else}). The proof of Theorem \ref{strong-c-posit-th} itself, as a whole, is provided in Section \ref{sec-synth}, by synthetizing the various, sometimes redundant, smallness assumptions made in the previous sections on $|K-1|_{C^2(S)},\ \varsigma$ and an extra parameter $\delta$ used to locate $\exp_m(V)$ with respect to $m$ and $m^*$ as just described. The proof of the main perturbation lemma is deferred to Appendix \ref{app-a}, but Section \ref{sec-perturb} includes a straightforward application of it to a uniform convexity estimate for the boundary of $\noconj$.\\

Finally, a warning must be made about some notations and conventions used below. Starting from Lemma \ref{mainperturb-lem} (Section \ref{sec-perturb}), we will abbreviate $\ds |K-1|_{C^2(S)}$ merely by $\varepsilon$. In Section \ref{sec-conj} (resp. Section \ref{sec-zero}), we will set $\delta_1 d(m,m^*)$ (resp. $\delta_2$) for the maximal distance assumed between $\exp_m(V)$ and the first conjugate point $m^*$  (resp. and the point $m$); consistently in Section \ref{sec-else}, we will set $\frac{1}{2}\delta_1 d(m,m^*)$ (resp. $\frac{1}{2}\delta_2$) for the minimal distance at which $\exp_m(V)$ must stay away from $m^*$ (resp. from $m$) on that geodesic. In the course of the proof, starting from Lemma \ref{mainperturb-lem}, we will require various (fairly explicit, universal) smallness conditions on $\varepsilon$ or the auxiliary position parameters $\delta_i$'s. Furthermore, in each case or subcase distinguished below for $(m,V,\xi,\nu)$, we will find a different value of the (small positive) constant $\varsigma$ occuring in (\ref{strong-c-posit-ineq}); the actual value to be taken for $\varsigma$ in the statement of Theorem \ref{strong-c-posit-th} will be, of course, the \emph{smallest} among them. The various universal\footnote{thus, in particular, independent of $(m,V)\in\noconj$} constants and smallness conditions arising in the paper are listed in Appendix \ref{app-b} to which the reader should systematically refer.

\section{$c$-curvature expression in dimension 2}\label{sec-express}
Henceforth, we fix a surface $S$, a point $m_0\in S$ and three non-zero tangent vectors
$(V_0,\xi,\nu)$ in
$T_{m_0}S$ with
$(m_0,V_0)\in\nocut$ and
$(\xi,\nu)$ linearly independent. We wish to calculate the $c$-curvature ${\cal C}(m_0,V_0)(\xi,\nu)$.
\subsection{General case}\label{sec-genexpress}
A chart $x=(x^1,x^2)$ of $S$ centered at
$m_0$ such that the local components
$g_{ij}(x)$ of the metric satisfy:
$g_{ij}(0)=\delta_{ij},\ dg_{ij}(0)=0$, is called \textit{normal} at $m_0$; let $x$ be such a chart. We set $v=(v^1,v^2)$ for the fiber
coordinates of
$TS\to S$ naturally associated to $x$, use Einstein's convention and abbreviate partial derivatives as follows:
$$\partial_i=\frac{\partial}{\partial x^i},\partial_{ij}=\frac{\partial^2}{\partial x^i\partial x^j},\ldots;D_i=\frac{\partial}{\partial
v^i},D_{ij}=\frac{\partial^2}{\partial v^i\partial v^j},\ldots$$
For each $(m,V)\in\nocut$ with $m$ in the domain of the chart $x$, we set:
$$X=X(x,v,t)=\left(X^1(x,v,t),X^2(x,v,t)\right)=x\left(\exp_m(tV)\right),$$
where $x=x(m)$ and $V=v^i\partial_i$. For $V\in T_{m_0}S$ such that $(m_0,V)\in\nocut$, and setting $\xi=\xi^i\partial_i$, we recall from
\cite{delge} that the quadratic form $A(m_0,V)(\xi)$ defined in the introduction is equal to
$A_{ij}(v)\xi^i\xi^j$ with:
\begin{equation}\label{hessian-c-eq}
A_{ij}(v)=Y^i_k(v)\ \partial_jX^k(0,v,1)
\end{equation}
and the matrix $Y^i_k(v)$ given by: $Y^i_k(v)\ D_jX^k(0,v,1)=\delta^i_j$. Given $V=v^i\partial_i$ as above, it is convenient to compute the
right-hand side of (\ref{hessian-c-eq}) by choosing for $x$ a particular normal chart at $m_0$ (unique up to $x^1\to -x^1$), namely:
\begin{Definition}\label{fermi-chart-def}
A Fermi chart along $V$ is a normal chart $x$ at $m_0$ such that $V=r\partial_2$ (with $r=|V|$) and the Riemannian metric reads:
$$g=dx^1\otimes dx^1+G(x^1,x^2)\ dx^2\otimes dx^2,\ {\it with}\ G(0,x^2)=1,\ \partial_1G(0,x^2)=0.$$
\end{Definition}
Let $x$ be a Fermi chart along $V$. The geodesic $t\in[0,1]\to m_t=\exp_m(tV)\in S$ (called the \textit{axis} of the chart) simply reads
$t\mapsto X((0,0),(0,r),t)=(0,tr)$ and, for fixed $x^2$, the paths which read $t\mapsto (t,x^2)$ are geodesics of $S$ as well, orthogonal to
the axis. The Christoffel symbols are given by:
$$\Gamma^1_{22}=-\frac{1}{2}\partial_1G,\ \Gamma^2_{12}=\frac{\partial_1G}{2G},\ \Gamma^2_{22}=\frac{\partial_2G}{2G},\ {\rm others\ vanish},$$
and the Gauss curvature, by $\ds K=-\frac{\partial_{11}\sqrt{G}}{\sqrt{G}}$. We thus get for the derivatives of the Christoffel
symbols on the axis, intrinsic expressions given in terms of $K$ at $x=(0,x^2)$ by:
$$\partial_1\Gamma^1_{22}=-\partial_1\Gamma^2_{12}=K,\ \partial_{11}\Gamma^1_{22}=-\partial_{11}\Gamma^2_{12}=\partial_1K,\
\partial_1\Gamma^2_{22}=0,\ \partial_{11}\Gamma^2_{22}=-\partial_2K.$$
With these formulas at hand, we readily find:
$$\partial X((0,0),(0,r),t)=\left( \begin{array}{cc}
f_0(t) & 0 \\
0 & 1
\end{array}\right),\ D X((0,0),(0,r),t)=\left( \begin{array}{cc}
f_1(t) & 0 \\
0 & t
\end{array}\right),$$
where $f_i(t)=f_i((0,0),(0,r),t)$ for $i\in\{0,1\}$; here, $f_i(x,w,t)$ are the expressions in the chart $x$ of the solutions for $t\in[0,1]$ of
the Jacobi equation:
\begin{equation}\label{jacob-eq}
\ddot{f}+|W|^2\ K\left(\exp_m(tW)\right) f=0
\end{equation}
(where $x=x(m), W=w^i\partial_i$ with $(m,W)\in\noconj$, and we use the dot notation: $\ds\dot{f}=\frac{df}{dt},\ \ddot{f}=\frac{d^2f}{dt^2}$),
satisfying the initial condition:
$$f_i(0)=\delta_{i0},\ \dot{f}_i(0)=\delta_{i1}\ .$$
\begin{Remark}\label{sturm-rk}
{\rm For later use, we observe that, for $t\in(0,1]$ and $(m,W)\in\noconj$, we have: $0<f_1(x,w,t)$. Moreover, Sturm
comparison theorem
\cite{carm} combined with Remark \ref{comparison-rk} provides the pinching:
$$\frac{\sin\left(\sqrt{\max_SK}|W|t\right)}{\sqrt{\max_SK}|W|}\leq f_1(x,w,t)\leq \frac{\sin(|W|t)}{|W|},$$
which yields $f_1(x,w,t)\leq t\leq 1$ and $\ds\lim_{|W|\downarrow 0}f_1(x,w,1)=1$.}
\end{Remark}
Back to $(m,W)=(m_0,V)$, applying
(\ref{hessian-c-eq}) in our Fermi chart along $V$, we get:
$$A(m_0,V)(\xi)=|\xi|^2-\left(1-\frac{f_0(1)}{f_1(1)}\right) |\xi-g(\xi,U)U|^2,\ {\rm with}\ U=\frac{V}{|V|}.$$
Here comes a key observation, also made in \cite{firivi1} (and extended to the higher dimensional setting in \cite{firivi2}, see also \cite{leemcc}): the right-hand
side of the preceding equation is intrinsic because so is (\ref{jacob-eq}). We may thus use a single Fermi chart $x$, along the
sole tangent vector
$V_0$ at
$m_0$, and write for each
$V=v^i\partial_i\in T_{m_0}S$ close to
$V_0$:
\begin{equation}\label{fermi-hess-c-eq}
A(m_0,V)(\xi)=|\xi|^2-\left(1-\frac{f_0(0,v,1)}{f_1(0,v,1)}\right) |\xi-g(\xi,U)U|^2.
\end{equation}
We will now calculate the $c$-curvature ${\cal C}(m_0,V_0)(\xi,\nu)$ in that Fermi chart (fixed once for all), by combining
(\ref{c-curvdef-eq}) with (\ref{fermi-hess-c-eq}). Letting henceforth $\xi$ and $\nu$ be unit vectors and orienting the tangent plane
$T_{m_0}S$ by the local volume form $dx^1\wedge dx^2$, we denote by
$\vartheta$ (resp.
$\varphi$) the angle in $[0,2\pi)$ by which a direct rotation brings $\xi$ (resp. $\nu$) to $\ds U_0=\frac{V_0}{|V_0|}=\partial_2$; in other
words, we set:
$$\xi=\sin\vartheta\ \partial_1+\cos\vartheta\ \partial_2,\ \nu=\sin\varphi\ \partial_1+\cos\varphi\ \partial_2.$$
A lengthy but routine calculation yields:
\begin{eqnarray}\label{c-curv-general-eq}
{\cal C}(m_0,V_0)(\xi,\nu) = &-& \sin^2\vartheta\left(
\frac{f_0''}{f_1}-\frac{f_0f_1''}{f_1^2}-\frac{2f_0'f_1'}{f_1^2}+\frac{2f_0(f_1')^2}{f_1^3}\right)\\ &+&
\frac{2}{r_0^2}\left(\cos^2\vartheta-\cos^2(\vartheta+\varphi)\right)\left(1-\frac{f_0}{f_1}\right)\nonumber\\ &+&
\frac{4}{r_0}\cos\vartheta\sin\vartheta\sin\varphi\left( \frac{f_0'}{f_1}-\frac{f_0f_1'}{f_1^2}\right),\nonumber
\end{eqnarray}
where we have set, for short: $f_a'=\nu^iD_if_a(0,v_0,1),\ f_a''=\nu^i\nu^jD_{ij}f_a(0,v_0,1),$
for
$a=0,1$, and
$v_0=(0,r_0)$ with
$r_0=|V_0|$.
%%%
\subsection{Constant curvature case recalled}\label{spherical-subsec}
Setting for short $\kappa=K(m_0)$ and $\bar{r}=\sqrt{\kappa}\ r$, let us recall the expressions which occur for $f_0, f_1$ in case
$K\equiv\kappa$, labelling them all with a bar:
$$\bar{f}_0(0,v,t)=\cos(\bar{r}t),\ \bar{f}_1(0,v,t)=\frac{\sin(\bar{r}t)}{\bar{r}},\ {\rm where}\ r=\sqrt{(v^1)^2+(v^2)^2}.$$
At $(v,t)=(v_0,1)$, with $v_0=(0,r_0)$ and $\bar{r}_0=\sqrt{\kappa}\ r_0$, we infer correspondingly:
$$\bar{f}_0'=-\sqrt{\kappa}\ \sin\bar{r}_0\ \cos\varphi,\ \
\bar{f}_0''=\kappa\left(-\frac{\sin\bar{r}_0}{\bar{r}_0}+\left(\frac{\sin\bar{r}_0}{\bar{r}_0}-\cos\bar{r_0}\right)\cos^2\varphi\right),$$
$$\bar{f}_1'=\frac{\sqrt{\kappa}}{\bar{r}_0}\left(\cos\bar{r}_0-\frac{\sin\bar{r}_0}{\bar{r}_0}\right)\cos\varphi\ ,$$
$$
\bar{f}_1'' =\frac{\kappa}{\bar{r}_0^2}\left(\cos\bar{r}_0-\frac{\sin\bar{r}_0}{\bar{r}_0}+\left(3\left(\frac{\sin\bar{r}_0}{\bar{r}_0}-\cos\bar{r}_0
\right)-\bar{r}_0\sin\bar{r}_0
\right)\cos^2\varphi\right),$$
hence:
\begin{eqnarray}\label{c-curv-sphere-eq}
\frac{1}{\kappa}\ \overline{\cal C}(m_0,V_0)(\xi,\nu) &=&
\sin^2\vartheta\ \sin^2\varphi\ \frac{\bar{r}_0^2+\bar{r}_0\cos\bar{r}_0\sin\bar{r}_0-2\sin^2\bar{r}_0}{\bar{r}_0^2\sin^2\bar{r}_0}\\
&+& 2 \sin^2\vartheta\ \cos^2\varphi\ \frac{\sin\bar{r}_0-\bar{r}_0\cos\bar{r}_0}{\sin^3\bar{r}_0}\nonumber\\
&+& 2\cos^2\vartheta\ \sin^2\varphi\
\frac{\sin\bar{r}_0-\bar{r}_0\cos\bar{r}_0}{\bar{r}_0^2\sin\bar{r}_0}\nonumber\\
&+& 4\cos\vartheta\ \sin\vartheta\ \cos\varphi³\ \sin\varphi\ \frac{\sin^2\bar{r}_0-\bar{r}_0^2}{\bar{r}_0^2\sin^2\bar{r}_0}.\nonumber
\end{eqnarray}
\section{Perturbative tools}\label{sec-perturb}
In the sequel of the paper, dropping the first argument $x=x(m)$ since it is fixed, equal to $(0,0)=x(m_0)$, we simply write:
$f_a=f_a(v,t),\ X=X(v,t)$ and, abusively with the same letter: $K(X(v,t))=K\left( \exp_{m_0}(tV)\right)$, where
$V=v^i\partial_i$. Moreover, anytime the second argument $v$ is equal to $v_0=(0,r_0)$, we will also drop it and just write: $f_a=f_a(t)$ and
so on.\\Given a real number
$\omega>0$, we will require the linear map:
$$f\in C^0([0,1],{\Bbb R})\longrightarrow{\cal S}_{\omega}(f)\in C^0([0,1],{\Bbb R})$$
defined as the solution map $f\mapsto u$ of the linear initial
value problem:
$$\ddot{u}+\omega^2u=f,\ u(0)=\dot{u}(0)=0.$$
The representation formula : $\ds {\cal S}_{\omega}(f)(t)=\int_0^t\frac{\sin\left(\omega(t-\tau)\right)}{\omega}\ f(\tau)\
d\tau$ is well known. Setting $\ds \Vert v\Vert=\sup_{t\in[0,1]}|v(t)|$, it yields for ${\cal S}_{\omega}$ the
contraction estimate:
\begin{equation}\label{contraction-estim}
\Vert {\cal S}_{\omega}(f) \Vert\leq \frac{1}{2}\Vert f\Vert,
\end{equation}
easily obtained by writing:
$$u(t)=\int_0^t\dot{u}(\tau)d\tau=\int_0^t\int_0^{\tau}\cos\left(\omega(\tau-\theta)\right)\ f(\theta)\ d\theta d\tau .$$
We will also require the following formulas (written at $t=1$, for
$f(t)=t$ and
$f(t)=t^2$):
\begin{equation}\label{auxil-sol-formulas}
{\cal S}_{\bar{r}_0}(t)(1)=\frac{\bar{r}_0-\sin\bar{r}_0}{\bar{r}_0^3},\ {\cal S}_{\bar{r}_0}(t^2)(1)=\frac{\bar{r}_0^2+2(\cos\bar{r}_0-1)}{\bar{r}_0^4}.
\end{equation}
We are now ready to state our main perturbation lemma, the proof of which is deferred to Appendix \ref{app-a}:
\begin{Lemma}\label{mainperturb-lem}
If
$\ds |K-1|_{C^2(S)}\leq\frac{1}{\pi^2}$, there exists universal constants $B_{1ka}, B_{2ka}, B_{3ka},$ for $a\in\{0,1\}$ and $k\in \{0,1,2\}$, such that the following estimates hold:
$$\Vert D_{\nu}^kf_a\Vert\leq B_{1ka},\ \Vert D_{\nu}^k(f_a-\bar{f}_a)\Vert\leq B_{2ka}\ \varepsilon\ r_0^{2-k},$$
$$\Vert D_{\nu}^k(f_a-\bar{f}_a)+r_0^{3-k}\ \psi_k\ {\cal S}_{\bar{r}_0}(t^{a+1})\Vert \leq B_{3ka}\ \varepsilon\ r_0^{4-k}\ ,$$
where, for short, $\varepsilon:=|K-1|_{C^2(S)}$ and:
$$\psi_0:=\partial_2K(0),\ \psi_1:=3\cos\varphi\ \partial_2K(0)+\sin\varphi\ \partial_1K(0),$$
$$\psi_2:= (2+4\cos^2\varphi)\ \partial_2K(0)+4\sin\varphi\cos\varphi\ \partial_1K(0)$$
(from now on, we will freely use to these abbreviations).
\end{Lemma}
\begin{Remark}\label{why-lem-rem}
{\rm Let us stress that the
bounds:
$$\forall a=0,1,\ \Vert D_{12}f_a\Vert\leq 2 B_{12a},\ \Vert D_{12}(f_a-\bar{f}_a)\Vert\leq 2 B_{22a}\varepsilon,$$
follow from thoses on $\Vert D_{\nu\nu}f_a\Vert$ and $\Vert D_{\nu\nu}(f_a-\bar{f}_a)\Vert$ by letting $\nu=\frac{1}{\sqrt{2}}(\partial_1+\partial_2)$.}
\end{Remark}

The first line of conclusion of Lemma \ref{mainperturb-lem} will be used to prove Theorem \ref{th-main} near\footnote{where $\overline{\cal C}(m_0,V_0)(\xi,\nu)$ could blow up since $\bar{r}_0$ could exit from $(0,\pi)$ for $\varepsilon\not=0$}
the first conjugate point (Section \ref{sec-conj}). Uniformly away from that point, and crucially for $r_0\downarrow 0$, the proof requires the second line of conclusion through a Maclaurin type approximation estimate for the $c$-curvature, namely:
\begin{Corollary}\label{c-curv-perturb-cor}
If
$\ds |K-1|_{C^2(S)}\leq\frac{1}{\pi^2}$ and $\bar{r}_0<\pi$, there exists a universal constant $C_1$ such that the absolute value of the following expression:
\begin{eqnarray}
\frac{f_1^3}{\bar{f}_1^3}\ {\cal C}(m_0,V_0)(\xi,\nu) &-& \overline{\cal C}(m_0,V_0)(\xi,\nu)- \frac{r_0\psi_2\sin^2\vartheta}{\bar{f_1}}\left( {\cal S}_{\bar{r}_0}(t)(1)-\frac{\bar{f}_0{\cal S}_{\bar{r}_0}(t^2)(1)}{\bar{f_1}}\right)\nonumber\\
&+& \frac{2r_0\psi_0{\cal S}_{\bar{r}_0}(t^2-t)(1)}{\bar{f_1}}\left(\cos^2\vartheta - \cos^2(\vartheta+\varphi)\right)\nonumber\\
&+& \frac{4r_0\psi_1\cos\vartheta\sin\vartheta\sin\varphi}{\bar{f_1}}\left( {\cal S}_{\bar{r}_0}(t)(1)-\frac{\bar{f}_0{\cal S}_{\bar{r}_0}(t^2)(1)}{\bar{f_1}}\right)\nonumber
\end{eqnarray}
is bounded above by:
$$\frac{1}{\bar{f}_1^3}C_1^3\pi^8\varepsilon r_0^2(338\sin^2\vartheta +268 \cos^2\vartheta\sin^2\varphi).$$
\end{Corollary}
\textbf{Proof of the corollary.} For each $a\in\{0,1\}$ and $k\in \{0,1,2\}$, we split $D_{\nu}^kf_a$ identically into three summands: $D_{\nu}^kf_a=S_1^{(k,a)}+S_2^{(k,a)}+S_3^{(k,a)}$ given by:
$$S_1^{(k,a)}=D_{\nu}^k\bar{f}_a,\ S_2^{(k,a)}=-r_0^{3-k}\psi_k{\cal S}_{\bar{r}_0}(t^{a+1})(1).$$
From (\ref{auxil-sol-formulas}), we define the constants $c_{6},c_{7}$ as in Appendix \ref{app-b}. From Lemma \ref{mainperturb-lem}, we know that
$$\left|S_1^{(k,a)}\right|\leq B_{1ka},\ \left|S_3^{(k,a)}\right|\leq B_{3ka}\varepsilon r_0^{4-k},$$
and from the obvious bounds:
\begin{equation}\label{obv-psi-bound}
\left|\psi_0\right|\leq \varepsilon,\ \left|\psi_1\right|\leq 4\varepsilon,\ \left|\psi_2\right|\leq 8\varepsilon,
\end{equation}
we further know that
$$\left|S_2^{(k,a)}\right|\leq 8c_{6+a}\varepsilon r_0^{3-k}\ .$$
Let us consider the expression (\ref{c-curv-general-eq}) of the $c$-curvature, multiply it by $f_1^3$ and, using the preceding splittings and bounds, let us estimate the Maclaurin approximation of each of the three auxiliary expressions:
$$E_1:=f_1^2f_0''-f_0f_1f_1''-2f_1f_0'f_1'+2f_0(f_1')^2$$
$$E_2:=\frac{2}{r_0^2}f_1^2\left(f_1-f_0\right),\ E_3:=\frac{4}{r_0}f_1\left( f_1f_0'-f_0f_1'\right),$$
which occur in $f_1^3{\cal C}(m_0,V_0)(\xi,\nu)$ as coefficients, respectively of:
$$-\sin^2\vartheta,\ \left(\cos^2\vartheta-\cos^2(\vartheta+\varphi)\right),\ \cos\vartheta\sin\vartheta\sin\varphi\ .$$
Setting $\overline{E}_1,\overline{E}_2,\overline{E}_3,$ for the corresponding quantities defined with $\bar{f}_0,\bar{f}_1$ instead of $f_0,f_1$, and proceeding stepwise, with careful intermediate calculations\footnote{in particular, for counting numbers of terms which are $O(\varepsilon r_0^2)$}, we get for the $(E_\ell-\overline{E}_\ell)$'s the following analogues of the second line of conclusion of Lemma \ref{mainperturb-lem}:
$$\left\vert E_1-\overline{E}_1+r_0\psi_2\bar{f}_1\left[\bar{f}_1 {\cal S}_{\bar{r}_0}(t)(1) - \bar{f}_0{\cal S}_{\bar{r}_0}(t^2)(1)\right]\right\vert\leq 154\pi^8C_1^3\varepsilon r_0^2\ ,$$
$$\left\vert E_2-\overline{E}_2+2r_0\psi_0\bar{f}_1^2 {\cal S}_{\bar{r}_0}(t^2-t)(1)\right\vert\leq 84\pi^8C_1^3\varepsilon r_0^2\ ,$$
$$\left\vert E_3-\overline{E}_3+4r_0\psi_1\bar{f}_1\left[\bar{f}_1 {\cal S}_{\bar{r}_0}(t)(1) - \bar{f}_0{\cal S}_{\bar{r}_0}(t^2)(1)\right]\right\vert\leq 200\pi^8C_1^3\varepsilon r_0^2\ ,$$
where the constant $C_1$ is defined\footnote{using, in particular, the bounds $\sqrt{\kappa}\leq 1+\frac{1}{2\pi^2}<\frac{19}{18}$ and $\kappa\leq 1+\frac{1}{\pi^2}<\frac{10}{9}$} in Appendix \ref{app-b}, as well as three other constants $c_8,c_9,c_{10}$, and where, recalling Remark \ref{comparison-rk}, $\pi^8$ is used as an upper bound for $\max\left(1,r_0^{p-2}\right)$ with\footnote{for instance, $p=2$ (resp. $p=10$) for $\left(S_1^{(0,1)}\right)^2S_3^{(2,0)}$ (resp. $\left(S_3^{(0,1)}\right)^2S_3^{(2,0)}$) in the first term of $E_1$} $2\leq p\leq 10$. Since $\bar{r}_0<\pi$, we may divide by $\bar{f}_1^3>0$ the resulting Maclaurin approximation estimate for $f_1^3{\cal C}(m_0,V_0)(\xi,\nu)$ and, using the general inequalities:
\begin{eqnarray}\label{auxil-trigo-bound}
\left| \cos^2\vartheta - \cos^2(\vartheta+\varphi)\right| &\leq& \sin^2\vartheta +2 \cos^2\vartheta\sin^2\varphi\ ,\\
2|\cos\vartheta\sin\vartheta\sin\varphi| &\leq& \sin^2\vartheta +\cos^2\vartheta\sin^2\varphi\ ,\nonumber
\end{eqnarray}
we obtain the estimate of Corollary \ref{c-curv-perturb-cor}.

\paragraph{Quick digression on the convexity of $\noconj$.}
The reader may wish to skip the rest of this section,
devoted to a quick digression from our main topic.
Indeed, let us pause and
provide a
uniform convexity estimate on the tangential domains
$$\noconj_{m}=\{W\in T_mS, (m,W)\in\noconj\},$$
obtained in terms of $|K-1|_{C^2(S)}$ as a direct consequence of Lemma \ref{mainperturb-lem}, and stated as follows:
\begin{Corollary}\label{convex-conj-cor}
Let $S$ be a surface as above with: $\min_S K=1$. There exists universal positive constants $\beta,\gamma,C,$ with
$\ds\beta\leq\frac{1}{\pi^2}$ and $\gamma\leq C$, such that, if
$|K-1|_{C^2(S)}\leq\beta$, for each $m_0\in S$ and $V_0\in\partial\noconj_{m_0}$, the curvature of the boundary curve
$\partial\noconj_{m_0}$ at $V_0$ is pinched between $\gamma$ and $C$.
\end{Corollary}

Qualitative proofs of the uniform convexity of $\noconj$ are given in \cite{cari,firivi2} for $C^4$ perturbations of the standard $n$-sphere. Let us further note that, combining Corollary \ref{convex-conj-cor} with Theorem \ref{th-main}, one can readily show that $\nocut$ is convex for small enough $\beta$ by arguing as in \cite{firivi2}, here just with a linear path $t\in[0,1]\to V_t=tV_1+(1-t)V_0$ in $T_{m_0}S$, with $V_0$ and $V_1$ in $\nocut_{m_0}$.\\

\noindent \textbf{Proof.} Fix $(m_0,V_0)$ as stated and take a Fermi chart $x$ along $V_0$, sticking to the above notations. From the
vanishing of $f_1(0,v_0,t)$ at $t=1$ combined with its positivity for $t\in(0,1)$ and the uniqueness of the solution of the initial (here
final, rather) value problem
\cite{pw}, we infer that $\dot{f}_1(0,v_0,1)<0$ hence also $D_2f_1(0,v_0,1)<0$, since $f_1(0,v_0,t)=f_1((0,0),(0,r_0),t)\equiv
tf_1((0,0),(0,tr_0),1)$. Therefore, near $V_0$, the curve $\partial\noconj_{m_0}$ admits the equation $v^2=h(v^1)$ with the function $h$
implicitly given by:
$$f_1\left((0,0),(v^1,h(v^1)),1\right)=0,\ {\rm and}\ h(0)=r_0.$$
Now, classically \cite{carm-surf}, the curvature \textbf{k} of $\partial\noconj_{m_0}$ at $V_0$ is equal to:
$${\rm \textbf{k}}=\frac{-h''(0)}{\left(1+h'(0)^2\right)^{3/2}}\equiv -
\frac{D_{11}f_1(D_2f_1)^2-2D_{12}f_1(D_1f_1)(D_2f_1)+D_{22}f_1(D_1f_1)^2}{\left((D_1f_1)^2+(D_2f_1)^2 \right)^{3/2}}.$$
Considering this formula, and since with $\varepsilon=0$ we would have $r_0=\pi, f_1=\bar{f}_1$ and ${\rm \textbf{k}}=\frac{1}{\pi^2}$,
the timeliness of Lemma \ref{mainperturb-lem} for our purpose is fully conceivable. For an effective proof, we first observe that, by Sturm theorem
\cite{carm}, $r_0$ is
pinched between
$\pi/\sqrt{1+\varepsilon}$ and
$\pi$; in particular, we have:
\begin{equation}\label{eps-pinch-conj}
\pi\left(1-\frac{\varepsilon}{2}\right)\leq \bar{r}_0\leq\pi\left(1+\frac{\varepsilon}{2}\right).
\end{equation}
At $(x,v,t)=(0,v_0,t)$, Lemma \ref{mainperturb-lem} and the formulas of Section \ref{spherical-subsec} imply:
$$D_2f_1\leq\frac{1}{r_0}\left(\cos\bar{r}_0 -\frac{\sin\bar{r}_0}{\bar{r}_0}\right)+B_{211}\varepsilon r_0\ ,$$
$$D_{11}f_1\leq\frac{1}{r_0^2}\left(\cos\bar{r}_0 -\frac{\sin\bar{r}_0}{\bar{r}_0}\right)+B_{221}\varepsilon\ ,$$
which, combined with the pinching of $r_0$
and standard bounds on the cosine and sine, yields:
$$D_2f_1\leq\frac{1}{r_0}\left(-1+\varepsilon\left(\frac{1}{2}+B_{211}\pi^2\right)+\frac{\varepsilon^2\pi^2}{8} \right),$$
$$D_{11}f_1\leq\frac{1}{r_0^2}\left(-1+\varepsilon\left(\frac{1}{2}+B_{221}\pi^2\right)+\frac{\varepsilon^2\pi^2}{8} \right).$$
So $\ds D_2f_1\leq -\frac{1}{2r_0}$, hence in particular $\ds|Df_1|\geq\frac{1}{2r_0}\geq\frac{1}{2\pi}$, and $\ds D_{11}f_1\leq
-\frac{1}{2r_0^2},$ provided
$\varepsilon$ is small enough. Moreover, still by Lemma
\ref{mainperturb-lem} and Section \ref{spherical-subsec}, we have at $(0,v_0,t)$:
$$|D_1f_1|\leq B_{211}\varepsilon r_0\leq B_{211}\pi\varepsilon\ .$$
The combination of Lemma \ref{mainperturb-lem} (including Remark \ref{why-lem-rem}) with the preceding bounds yields, on the one hand:
$${\rm \textbf{k}}\leq (2\pi)^36B_{121}B_{111}^2=:C ,$$
on the other hand:
$${\rm \textbf{k}}\geq\frac{1}{2B_{111}^3\sqrt{2}}\left( \frac{1}{8r_0^4}-\varepsilon\pi
B_{211}B_{121}\left(4B_{111}+\varepsilon\pi B_{211}\right)\right)$$ so that $\ds{\rm
\textbf{k}}\geq\gamma:=\frac{1}{32\pi^4\sqrt{2}B_{111}^3}$ for
$\varepsilon$ small enough. Altogether, this pinching of \textbf{k} holds provided we require $\varepsilon\leq\beta$ with $\beta$ the smallest
among the positive roots of the quadratic equations\footnote{which turns out to be that of the second equation}:
$$\frac{\pi^2}{8}\beta^2+\left(\frac{1}{2}+B_{221}\pi^2\right)\beta-\frac{1}{2}=0\ ,$$
$$\pi^2B_{121}B_{211}^2\beta^2+4\pi B_{111}B_{121}B_{211}\beta-\frac{1}{16\pi^4}=0\ .$$
Finally, a tedious but routine evaluation shows that we may take:
$$\gamma = 1.1\times 10^{-7},\ C=7.4\times 10^{7},\ \beta=4.5\times 10^{-10},$$
in the statement of Corollary \ref{convex-conj-cor}.
\section{$c$-curvature almost-positivity near conjugacy}\label{sec-conj}
In this section, we prove Theorem \ref{strong-c-posit-th} at $(m_0,V_0)\in\noconj$ and $(\xi,\nu)$ unit vectors of $T_{m_0}S$ in case the point
$\exp_{m_0}(V_0)$ is close to the first conjugate point $m^*_0$ of $m_0$ along the geodesic $t\in{\Bbb R}^+\to\exp_{m_0}(tV_0)\in S$. Specifically, setting $\l_0$ for the length of that geodesic curve from $m_0$ up to $m_0^*$, we
establish the following proposition:
\begin{Proposition}\label{posit-conj-prop}
There exists a triple of small (strictly) positive real numbers $(\eta_1,\delta_1,\varsigma_1)$ such that ${\cal C}(m_0,V_0)(\xi,\nu)$ satisfies the lower bound (\ref{strong-c-posit-ineq}) with $\varsigma=\varsigma_1$, provided $\varepsilon=|K-1|_{C^2(S)}\leq\eta_1$ and $(1-\delta_1)\l_0\leq |V_0|< \l_0$.
\end{Proposition}

\noindent\textbf{Proof.} Sticking to previous notations and recalling (\ref{eps-pinch-conj}), we infer from the pinching of $|V_0|$ the following ones (dropping the subscript of $\delta_1$):
\begin{equation}\label{pinch-r-rbar-conj}
\left(1-\frac{\varepsilon}{2}-\delta\right)\pi\leq r_0 \leq \pi\  {\rm and}\ \left(1-\frac{\varepsilon}{2}-\delta\right)\pi\leq \bar{r}_0 < \left(1+\frac{\varepsilon}{2}\right)\pi\ .
\end{equation}
We will assume:
\begin{equation}\label{grosspinch-r-rbar-conj}
r_0\in \left(\frac{5\pi}{6} , \pi\right)\ {\rm and}\ \bar{r}_0\in\left(\frac{5\pi}{6} , \frac{7\pi}{6}\right)
\end{equation}
with no loss of generality (it holds under the smallness condition posed below\footnote{we will say, for short, that they are (\ref{epsdel-ineq2})-small} on $\varepsilon$ and $\delta$, see (\ref{epsdel-ineq2}) of Appendix \ref{app-b}). Combining (\ref{pinch-r-rbar-conj}) with Remark \ref{sturm-rk}, the formulas of Section \ref{spherical-subsec} and the first line of conclusion of Lemma \ref{mainperturb-lem}, we derive the following set of inequalities:
\begin{equation}\label{00-pinch}
-1-\varepsilon \pi^2 B_{200}\leq f_0\leq -1+\varepsilon \pi^2 B_{200}+\left(\frac{\varepsilon}{2}+\delta\right)^2\frac{\pi^2}{2}
\end{equation}
\begin{equation}\label{01-pinch}
0<f_1\leq \frac{\sin r_0}{r_0}\leq \frac{\frac{\varepsilon}{2}+\delta}{1-\left(\frac{\varepsilon}{2}+\delta\right)}
\end{equation}
\begin{equation}\label{10-bound}
|f'_0|\leq \left(\frac{\varepsilon}{2}+\delta\right)\left(1+\frac{\varepsilon}{2}\right)\pi + \varepsilon\pi B_{210}.
\end{equation}
Furthermore, we derive two important lower bounds, namely:
\begin{Lemma}\label{posit-minor-lem1}
If $\cos\varphi\not=0$, and $\varepsilon$ and $\delta$ satisfy the relative smallness condition\footnote{to be used only in Section \ref{subsec-cosnot0} below, with $|\cos\varphi|$ bounded away from 0 by a (small) universal constant \textit{i.e.} with $|\cos\varphi|$ \emph{replaced} by that constant} (\ref{epsdel-ineq1}) (see Appendix \ref{app-b}), the lower bound:
$$|f'_1|\geq \frac{|\cos\varphi|}{\pi} - \varepsilon\pi B_{211}-\frac{\frac{\varepsilon}{2}}{\pi\left(1-\left(\frac{\varepsilon}{2}+\delta \right)\right)}-\left(\frac{\varepsilon}{2}+\delta \right)^2\frac{\pi}{2}$$
holds true, as well as the sign condition:
$-f'_1\cos\varphi >0$ .
If $\ds |\cos\varphi|\leq\frac{1}{2}$ and $\varepsilon$ and $\delta$ are (\ref{epsdel-ineq2})-small, the following lower bound is valid instead:
$$f_0f_1''\geq\frac{1}{8\pi^2}-\varepsilon\left(B_{221}+\frac{1}{8}B_{200}\right)-\frac{1}{16}\left(\frac{\varepsilon}{2}+\delta\right)^2.$$
\end{Lemma}
\textbf{Proof.} If $\cos\varphi\not=0$, using $|f'_1-\bar{f}'_1|\leq\varepsilon\pi B_{211}$ combined with the lower bound:
$$-\frac{\bar{f}'_1}{\cos\varphi} \geq \frac{1}{\pi}\left(1-\frac{\frac{\varepsilon}{2}}{1-\left(\frac{\varepsilon}{2}+\delta\right)}-\frac{\pi^2}{2}\left(\frac{\varepsilon}{2}+\delta\right)^2\right),$$
one can readily check the first part of the lemma. For the second part, we first note that $\bar{f}_1''$ is bounded above by the expression:
$$\frac{1}{r_0^2}\left[ -1+ \frac{\frac{\varepsilon}{2}}{1-\left(\frac{\varepsilon}{2}+\delta\right)}+\frac{\pi^2}{2}\left(\frac{\varepsilon}{2}+\delta\right)^2+\cos^2\varphi\left(\frac{3}{1-\left(\frac{\varepsilon}{2}+\delta\right)} + \frac{\varepsilon}{2}\left(1+\frac{\varepsilon}{2}\right)\pi^2\right)\right].$$
If $\ds |\cos\varphi| \leq \frac{1}{2}$,
it implies $\ds \bar{f}_1'' \leq -\frac{1}{8\pi^2}$ provided $\varepsilon$ and $\delta$ are taken (\ref{epsdel-ineq2})-small. By Lemma \ref{mainperturb-lem}, the inequality
$$f_1''\leq -\frac{1}{8\pi^2}+\varepsilon B_{221}$$
follows. Combined with (\ref{00-pinch}), it yields the second part of the lemma.\medskip

In order to investigate the sign of the $c$-curvature expression (\ref{c-curv-general-eq}), we will have to recast this expression in appropriate forms, namely, either:\smallskip
\begin{eqnarray}\label{altA-c-curv-eq}
{\cal C}(m_0,V_0)(\xi,\nu) = &-& \sin^2\vartheta\left(
\frac{f_0''}{f_1}-\frac{f_0f_1''}{f_1^2}-\frac{2f_0'f_1'}{f_1^2}\right)\nonumber\\ &-&
\frac{2f_0}{f_1}\left( \frac{f'_1}{f_1}\sin\vartheta+\frac{\sin\varphi\cos\vartheta}{r_0}\right)^2\\
&+&
\frac{2}{r_0^2}\left(1-\frac{f_0}{f_1}\right)\left(2\cos\vartheta\cos\varphi\sin\vartheta\sin\varphi-\sin^2\vartheta\sin^2\varphi\right)\nonumber\\ &+&
\frac{2}{r_0^2}\cos^2\vartheta\sin^2\varphi +\frac{4}{r_0}\cos\vartheta\sin\vartheta\sin\varphi\frac{f_0'}{f_1}\ ,\nonumber
\end{eqnarray}
or:
\begin{eqnarray}\label{altB-c-curv-eq}
{\cal C}(m_0,V_0)(\xi,\nu) = &-& \sin^2\vartheta\left(
\frac{f_0''}{f_1}-\frac{f_0f_1''}{f_1^2}-\frac{2f_0'f_1'}{f_1^2}+\frac{2f_0(f_1')^2}{f_1^3}+\frac{2}{r_0^2}\left(1-\frac{f_0}{f_1}\right)\right)\nonumber\\ &+&
\frac{2}{r_0^2}\left(1-\frac{f_0}{f_1}\right)\left(\cos\varphi\sin\vartheta+\cos\vartheta\sin\varphi\right)^2\\ &+&
\frac{4}{r_0}\cos\vartheta\sin\vartheta\sin\varphi\left( \frac{f_0'}{f_1}-\frac{f_0f_1'}{f_1^2}\right),\nonumber
\end{eqnarray}

\noindent We will also have to distinguish cases, depending on the size of $|\cos\varphi |$, then on the relative size of further arising quantities. In each case, relying on Lemma \ref{posit-minor-lem1} and treating $f_1$ as a small parameter in intermediate steps thanks to (\ref{01-pinch}), we will be able to find a leading term blowing up positively as $\varepsilon$ and $\delta$  go to zero and argue with it.

We are now ready to continue the proof of Proposition \ref{posit-conj-prop} and start out for a case by case discussion of  the sign of the $c$-curvature.
\subsection{Case $|\cos\varphi|$ small enough}
\subsubsection{Subcase $\ds \left\vert \sin\vartheta\frac{f'_1}{f_1}\right\vert\leq\frac{|\sin\varphi\cos\vartheta |}{2r_0}$}
In this subcase, the assumption $\ds |\cos\varphi|\leq\frac{1}{2}$ will suffice. We note the estimate:
$$\left( \frac{f'_1}{f_1}\sin\vartheta+\frac{\sin\varphi\cos\vartheta}{r_0}\right)^2\geq \frac{\sin^2\varphi\cos^2\vartheta}{4r_0^4}$$
and use it to derive from (\ref{altA-c-curv-eq}) the inequality:
\begin{eqnarray}
{\cal C}(m_0,V_0)(\xi,\nu) & \geq &  \frac{\sin^2\vartheta}{f_1^2}\left[f_0f_1''+2f_0'f_1'-f_1\left( f_0'' +\frac{2}{r_0^2} (f_1-f_0)\sin^2\varphi\right)\right]\nonumber\\ &-&
\frac{2f_0}{f_1}\frac{\sin^2\varphi\cos^2\vartheta}{4r_0^2}\nonumber\\
&-&
\frac{4}{f_1}|\cos\vartheta\sin\vartheta\sin\varphi| \left( \frac{f_1-f_0}{r_0^2} + \frac{|f_0'|}{r_0} \right)\ .\nonumber
\end{eqnarray}
The right-hand side will be handled relying on the second part of Lemma \ref{posit-minor-lem1} combined with the pinching (\ref{grosspinch-r-rbar-conj}) of $r_0$ and previous estimates on the various $D^kf_a$ terms which arise apart from $f_0f''_1$. Doing so, we can establish for ${\cal C}(m_0,V_0)(\xi,\nu)$ the lower bound:
\begin{eqnarray}\label{c-curv-min-conj1}
{\cal C}(m_0,V_0)(\xi,\nu) & \geq &  \frac{\sin^2\vartheta}{f_1^2}\left(\frac{1}{12\pi^2} -R_1(\varepsilon,\delta)\right)+ \frac{\sin^2\varphi\cos^2\vartheta}{8\pi^2f_1}\\ &+&
\frac{\sin^2\vartheta}{24\pi^2f_1^2}+ \frac{\sin^2\varphi\cos^2\vartheta}{8\pi^2f_1}-4\left( 1+\frac{5}{\pi^2}\right)
\frac{|\cos\vartheta\sin\vartheta\sin\varphi|}{f_1}\ ,\nonumber
\end{eqnarray}
provided $\varepsilon$ and $\delta$ are (\ref{epsdel-ineq3})-small, where $R_1(\varepsilon,\delta)$ stands for the rational function of $(\varepsilon,\delta)$ vanishing at $(0,0)$ given by the right-hand side of the smallness condition (\ref{epsdel-ineq3}). We claim that the second line of the right-hand side of (\ref{c-curv-min-conj1}) is non-negative for small enough $\varepsilon$ and $\delta$. Indeed, from the identity $a^2+b^2\geq 2|ab|$ used with $\ds a=\frac{\sin\vartheta}{2\pi\sqrt{6} f_1}$ and $\ds b=\frac{\sin\varphi\cos\vartheta}{2\pi \sqrt{2f_1}}$, we infer that this line is bounded below by:
$$\frac{|\cos\vartheta\sin\vartheta\sin\varphi|}{f_1\sqrt{f_1}}\left( \frac{1}{4\sqrt{3}\pi^2}-4\left( 1+\frac{5}{\pi^2}\right)\sqrt{f_1}\right),$$
and the claim follows by taking $\varepsilon$ and $\delta$ (\ref{epsdel-ineq3bis})-small. Eventually, for $\varepsilon$ and $\delta$ (\ref{epsdel-ineq3})(\ref{epsdel-ineq3bis})-small, we obtain: $${\cal C}(m_0,V_0)(\xi,\nu)\geq \frac{\sin^2\vartheta}{24\pi^2f_1^2}+ \frac{\sin^2\varphi\cos^2\vartheta}{8\pi^2f_1}\ .$$
Combining this lower bound with (\ref{01-pinch})(\ref{epsdel-ineq3bis}) and the useful, easily established\footnote{hint: use Remark \ref{comparison-rk}}, inequality:
\begin{equation}\label{auxil-trigo-bound2}
\sin^2\vartheta + \cos^2\vartheta\sin^2\varphi\geq \frac{1}{4\pi^2}{\cal A}_2(m_0,V_0,\xi,\nu),
\end{equation}
we get (\ref{strong-c-posit-ineq}) at $(m_0,V_0,\xi,\nu)$ with $\ds\varsigma=18$.
\subsubsection{Subcase $\ds \left\vert \sin\vartheta\frac{f'_1}{f_1}\right\vert >\frac{|\sin\varphi\cos\vartheta |}{2r_0}$}
The second line of the right-hand side of (\ref{altA-c-curv-eq}) is non-negative due to (\ref{00-pinch})(\ref{01-pinch}). So we may write:
\begin{eqnarray}
{\cal C}(m_0,V_0)(\xi,\nu) & \geq &  \frac{\sin^2\vartheta}{f_1^2}\left[f_0f_1''+2f_0'f_1'-f_1\left( f_0'' +\frac{2}{r_0^2} (f_1-f_0)\sin^2\varphi\right)\right]\nonumber\\ &-&
\frac{4}{f_1}|\cos\vartheta\sin\vartheta\sin\varphi| \left( \frac{f_1-f_0}{r_0^2}|\cos\varphi| + \frac{|f_0'|}{r_0} \right)\ ,\nonumber
\end{eqnarray}
hence also:
\begin{eqnarray}
{\cal C}(m_0,V_0)(\xi,\nu) & \geq &  \frac{\sin^2\vartheta}{f_1^2}\left[f_0f_1''+2f_0'f_1'-f_1\left( f_0'' +\frac{2}{r_0^2} (f_1-f_0)\sin^2\varphi\right)\right]\nonumber\\ &-&
\frac{8\sin^2\vartheta}{f_1^2}|f'_1| \left( \frac{f_1-f_0}{r_0}|\cos\varphi| + |f_0'| \right)\nonumber
\end{eqnarray}
by applying our subcase assumption. If $|\cos\varphi| \leq \frac{1}{2}$,  repeating the above argument, we see that the first line of the preceding right-hand side is larger than $\ds \frac{\sin^2\vartheta}{f_1^2}\left(\frac{1}{8\pi^2}-R_1(\varepsilon,\delta)\right)$, while the second line is bounded below by:
$$-8B_{111}\frac{\sin^2\vartheta}{f_1^2}\left[ \left(\frac{\varepsilon}{2}+\delta\right)\left(1+\frac{\varepsilon}{2}\right)\pi + \varepsilon\pi B_{210} + \frac{ 1+\varepsilon\pi^2B_{200}+\frac{\frac{\varepsilon}{2}+\delta}{1-\left(\frac{\varepsilon}{2}+\delta\right)}}{\pi \left( 1-\frac{\varepsilon}{2}-\delta\right)}|\cos\varphi|\right]$$
as shown by combining Lemma \ref{mainperturb-lem} with (\ref{pinch-r-rbar-conj})(\ref{00-pinch})(\ref{01-pinch})(\ref{10-bound}). Altogether, we may write:
\begin{eqnarray}
{\cal C}(m_0,V_0)(\xi,\nu) & \geq &  \frac{\sin^2\vartheta}{f_1^2}\left( \frac{1}{12\pi^2}-R_2(\varepsilon,\delta)\right)\nonumber\\ &+&
\frac{\sin^2\vartheta}{f_1^2}\left(\frac{1}{24\pi^2}- 8B_{111} \frac{ 1+\varepsilon\pi^2B_{200}+\frac{\frac{\varepsilon}{2}+\delta}{1-\left(\frac{\varepsilon}{2}+\delta\right)}}{\pi \left( 1-\frac{\varepsilon}{2}-\delta\right)} |\cos\varphi| \right)\nonumber
\end{eqnarray}
with $R_2(\varepsilon,\delta)$ given by the right-hand side of (\ref{epsdel-ineq4}). We get from (\ref{epsdel-ineq3bis}):
$$ \frac{ 1+\varepsilon\pi^2B_{200}+\frac{\frac{\varepsilon}{2}+\delta}{1-\left(\frac{\varepsilon}{2}+\delta\right)}}{\left( 1-\frac{\varepsilon}{2}-\delta\right)} \leq\frac{768(\pi^2+5)^2+2\pi^2+1}{768(\pi^2+5)^2-1}<1.00013\ ;$$
besides, we have: $\pi B_{111}=5+\pi\sqrt{2}+3\pi^2\simeq 39,05 <40$. So the smallness conditions:
\begin{equation}\label{small-cosphi-cond}
|\cos\varphi|\leq \frac{1}{7704}
\end{equation}
and (\ref{epsdel-ineq4}) imply that $\ds{\cal C}(m_0,V_0)(\xi,\nu)\geq \frac{\sin^2\vartheta}{24\pi^2f_1^2}$. In our present subcase, the latter inequality yields:
$${\cal C}(m_0,V_0)(\xi,\nu)\geq \frac{\sin^2\vartheta}{48\pi^2f_1^2}+\frac{\cos^2\vartheta\sin^2\varphi}{192\pi^2r_0^2f_1'^2}\ .$$
On the one hand, from (\ref{01-pinch}) combined with (\ref{epsdel-ineq3bis}), we get $48\pi^2f_1^2\leq 2\times 10^{-8}$. On the other hand, combining Lemma \ref{mainperturb-lem} with (\ref{pinch-r-rbar-conj}) and (\ref{small-cosphi-cond}), we have: $\ds r_0|f_1'|\leq \varepsilon \pi^2 B_{211}+\frac{1}{3852}$. So we can arrange to have $192\pi^2r_0^2f_1'^2\leq 1$ by taking $\varepsilon$ (\ref{epsdel-ineq4bis})-small. Altogether, we may write
$${\cal C}(m_0,V_0)(\xi,\nu)\geq \sin^2\vartheta+\cos^2\vartheta\sin^2\varphi$$
and, from (\ref{auxil-trigo-bound2}), conclude that (\ref{strong-c-posit-ineq}) holds at $(m_0,V_0,\xi,\nu)$, indeed, with $\ds\varsigma=\frac{1}{4\pi^2}$.
\subsection{Case $|\cos\varphi|>\frac{1}{7704}$}\label{subsec-cosnot0}
In this case, the first part of Lemma \ref{posit-minor-lem1} implies:
\begin{equation}\label{11-posit-minor}
|f'_1|\geq\frac{1}{15408\pi}\ {\rm with}\ -f'_1\cos\varphi >0\ ,
\end{equation}
provided $\varepsilon$ and $\delta$ are (\ref{epsdel-ineq5})-small. Furthermore, if the latter are (\ref{epsdel-ineq2})(\ref{epsdel-ineq3})-small, we infer from (\ref{00-pinch}) the pinching:
\begin{equation}\label{gross-00-pinch}
\frac{1}{2}\leq -f_0 \leq \frac{3}{2}\ .
\end{equation}
which will be used repeatedly. 
\subsubsection{Subcase $\cos\vartheta\cos\varphi\sin\vartheta\sin\varphi\leq 0$}
Working with the expression (\ref{altB-c-curv-eq}) of ${\cal C}(m_0,V_0)(\xi,\nu)$, the second line of which is non-negative, and combining (\ref{11-posit-minor}) with (\ref{gross-00-pinch}), (\ref{grosspinch-r-rbar-conj}) and Lemma \ref{mainperturb-lem}, we get the inequality:
\begin{eqnarray}
{\cal C}(m_0,V_0)(\xi,\nu) &\geq& \frac{\sin^2\vartheta}{f_1^3}\left(\frac{1}{15408^2\pi^2}-
f_1\left(B_{120}+\frac{3}{2}B_{121}+2B_{110}B_{111}+\frac{36}{5\pi^2}\right)\right)\nonumber\\ &+&
\frac{2}{r_0^2}\left(1-\frac{f_0}{f_1}\right)\left(\cos\varphi\sin\vartheta+\cos\vartheta\sin\varphi\right)^2\nonumber\\ &+&
\frac{2}{r_0f_1^2}|\cos\vartheta\sin\vartheta\sin\varphi|\left( \frac{1}{15408\pi} -2f_1B_{110}\right).\nonumber
\end{eqnarray}
Recalling (\ref{01-pinch}) (\ref{gross-00-pinch}) and assuming that $\varepsilon$ and $\delta$ are (\ref{epsdel-ineq6})-small, we infer the lower bound:
\begin{eqnarray}
{\cal C}(m_0,V_0)(\xi,\nu) &\geq& \frac{\sin^2\vartheta}{2\pi^2 15408^2f_1^3}+\frac{2}{r_0^2}\left(1+\frac{1}{2f_1}\right)\cos^2\vartheta\sin^2\varphi\nonumber\\
&+& \frac{1}{15408\pi r_0f_1^2}|\cos\vartheta\sin\vartheta\sin\varphi|\left(1-\frac{4}{r_0}15408\pi f_1(f_1-f_0)\right)\nonumber,
\end{eqnarray}
the second line of the right-hand side of which is non-negative, as checked by combining Remark \ref{sturm-rk} with (\ref{pinch-r-rbar-conj}) (\ref{01-pinch}) (\ref{gross-00-pinch}) and (\ref{epsdel-ineq3bis}). Using (\ref{pinch-r-rbar-conj})(\ref{01-pinch})(\ref{epsdel-ineq3bis}) to treat its first line, we obtain the inequality
$${\cal C}(m_0,V_0)(\xi,\nu) \geq 17205 \left(\sin^2\vartheta+\cos^2\vartheta\sin^2\varphi\right)$$
which, recalling (\ref{auxil-trigo-bound2}), implies (\ref{strong-c-posit-ineq}) at $(m_0,V_0,\xi,\nu)$ with $\varsigma=435$.
\subsubsection{Subcase $\cos\vartheta\cos\varphi\sin\vartheta\sin\varphi > 0$}
Here, since $-f'_1\cos\varphi >0$, we know that the expressions $\ds\sin\vartheta\frac{f'_1}{f_1}$ and $\ds\frac{\sin\varphi\cos\vartheta}{r_0}$ have \emph{opposite} signs.
\paragraph{Case $\ds \left\vert \sin\vartheta\frac{f'_1}{f_1}\right\vert\leq\frac{4|\sin\varphi\cos\vartheta |}{5r_0}$ or $\ds \left\vert \sin\vartheta\frac{f'_1}{f_1}\right\vert\geq\frac{5|\sin\varphi\cos\vartheta |}{4r_0}$:}
If $a$ and $b$ are two real numbers such that:
$ab<0$ and $\ds |a|\leq\frac{4}{5}|b|$ or $\ds |b|\leq\frac{4}{5}|a|$, one can readily verify that they satisfy: $\ds (a+b)^2\geq \frac{1}{50}(a^2+b^2)$.
Using the expression (\ref{altA-c-curv-eq}) of ${\cal C}(m_0,V_0)(\xi,\nu)$, we apply the preceding estimate with $\ds a= \sin\vartheta\frac{f'_1}{f_1}, b=\frac{\sin\varphi\cos\vartheta}{r_0},$ and find the $c$-curvature bounded below by:
$$\sin^2\vartheta\left(\frac{-f_0 (f'_1)^2}{25f_1^3} - \frac{f_0''}{f_1}+\frac{f_0f_1''+2f_0'f_1'}{f_1^2}-\frac{2}{r_0^2}\sin^2\varphi\left(1-\frac{f_0}{f_1}\right)\right)$$
$$+\frac{4|\cos\vartheta\sin\vartheta\sin\varphi|}{f_1}\left( -\frac{f_0}{7704r_0^2} - \frac{|f'_0|}{r_0}\right)-\frac{f_0}{25f_1r_0^2}\cos^2\vartheta\sin^2\varphi\ ,$$
hence also, combining Lemma \ref{mainperturb-lem} with (\ref{grosspinch-r-rbar-conj})(\ref{11-posit-minor})(\ref{gross-00-pinch}) and (\ref{10-bound}), by:
{\small $$\frac{\sin^2\vartheta}{f_1^3}\left(\frac{1}{50\pi^2 15408^2}-
f_1\left(B_{120}+\frac{3}{2}B_{121}+2B_{110}B_{111}+\frac{36}{5\pi^2}\right)\right)$$}
{\small $$+ \frac{4|\cos\vartheta\sin\vartheta\sin\varphi|}{f_1}\left( \frac{1}{15408\pi^2} - \frac{6}{5}\left( 
\left(\frac{\varepsilon}{2}+\delta\right)\left(1+\frac{\varepsilon}{2}\right) + \varepsilon B_{210} \right)\right)+\frac{1}{50\pi^2f_1}\cos^2\vartheta\sin^2\varphi.$$}
Recalling (\ref{01-pinch}), we infer that:
$${\cal C}(m_0,V_0)(\xi,\nu)\geq \frac{\sin^2\vartheta}{100\pi^2 15408^2f_1^3}+ \frac{1}{50\pi^2f_1}\cos^2\vartheta\sin^2\varphi\ ,$$
provided $\varepsilon$ and $\delta$ are (\ref{epsdel-ineq7})-small. Recalling (\ref{01-pinch})(\ref{epsdel-ineq3bis}) and (\ref{auxil-trigo-bound2}), it yields (\ref{strong-c-posit-ineq}) with $\varsigma=8.7$ at $(m_0,V_0,\xi,\nu)$.
\paragraph{Case $\ds \frac{4|\sin\varphi\cos\vartheta |}{5r_0}<\left\vert \sin\vartheta\frac{f'_1}{f_1}\right\vert < \frac{5|\sin\varphi\cos\vartheta |}{4r_0}$:}
This case is more difficult because we cannot use the square occuring in the second line of (\ref{altA-c-curv-eq}) any more; all we can do now from (\ref{altA-c-curv-eq}) is write:
\begin{eqnarray}
{\cal C}(m_0,V_0)(\xi,\nu) &\geq& \sin^2\vartheta\left(- \frac{f_0''}{f_1}+\frac{f_0f_1''+2f_0'f_1'}{f_1^2}-\frac{2}{r_0^2}\sin^2\varphi\left(1-\frac{f_0}{f_1}\right)\right)\nonumber\\
&+& \frac{4}{r_0^2} |\sin\varphi\cos\vartheta | |\cos\varphi \sin\vartheta | \left( 1-\frac{f_0}{f_1}\right)\nonumber\\
&-& \frac{4}{r_0} |\sin\varphi\cos\vartheta | |\sin\vartheta | \frac{|f'_0|}{f_1}\nonumber
\end{eqnarray}
and, from our present assumption, infer for ${\cal C}(m_0,V_0)(\xi,\nu)$ the lower bound:
\begin{eqnarray}\label{c-curv-min-conj2}
&\sin^2\vartheta& \left(- \frac{f_0''}{f_1}+\frac{f_0f_1''+2f_0'f_1'}{f_1^2} - \frac{2}{r_0^2}\sin^2\varphi\left(1-\frac{f_0}{f_1}\right) \right)\nonumber\\
&-& \frac{16f_0}{5r_0f_1^2} \sin^2\vartheta  |f'_1\cos\varphi|- 5 \sin^2\vartheta \frac{|f'_1f'_0|}{f_1^2}\ .
\end{eqnarray}
We will factorize $\ds\frac{\sin^2\vartheta}{f_1^2}$ as leading blowing up term in this expression and seek a positive coefficient for it. Doing so, we focus on the terms:
$$-\frac{f_0\sin^2\vartheta}{f_1^2}\left( -f''_1+\frac{16}{5r_0}|f'_1\cos\varphi| \right),$$
thus carefully investigate the sign of the latter parenthesis. Using Lemma \ref{mainperturb-lem}, we find it bounded below by:
$$\left( -\bar{f}''_1+\frac{16}{5r_0}|\bar{f}'_1\cos\varphi| \right)-\varepsilon\left( B_{221}+\frac{16}{5}B_{211}\right).$$
Now, a direct calculation of $\ds \left( -\bar{f}''_1+\frac{16}{5r_0}|\bar{f}'_1\cos\varphi|\right)$, using the expressions of $\bar{f}'_1$ and $\bar{f}''_1$ given in Section \ref{spherical-subsec}, shows that it is equal to:
$$\frac{1}{r_0^2}\left[|\cos\bar{r}_0| \left(1+\frac{1}{5}\cos^2\varphi\right)+\frac{\sin\bar{r}_0}{\bar{r}_0}\left(1+\left(\bar{r}^2_0+\frac{1}{5}\right)\cos^2\varphi\right)\right];$$
recalling (\ref{pinch-r-rbar-conj}), we see that it will meet the required positivity. Back to the lower bound (\ref{c-curv-min-conj2}), rewritten as $\ds{\cal C}(m_0,V_0)(\xi,\nu) \geq \frac{\sin^2\vartheta}{f_1^2}E$ with $E$ equal to:
$$(-f_0)\left( -f''_1+\frac{16}{5r_0}|f'_1\cos\varphi| \right)-7|f'_0f'_1|-f_1\left( |f''_0|+\frac{2\sin^2\varphi}{r_0^2}(f_1-f_0)\right),$$
the preceding argument, combined with Lemma \ref{mainperturb-lem}, Remark \ref{sturm-rk} and (\ref{pinch-r-rbar-conj}) (\ref{grosspinch-r-rbar-conj}) (\ref{01-pinch}) (\ref{10-bound}) (\ref{gross-00-pinch}), implies that $\ds{\cal C}(m_0,V_0)(\xi,\nu)\geq\frac{\sqrt{3}\sin^2\vartheta}{8\pi^2f_1^2}$ provided $\varepsilon$ and $\delta$ are (\ref{epsdel-ineq8})-small. In the present subcase, the latter inequality implies:
$${\cal C}(m_0,V_0)(\xi,\nu)\geq\frac{\sqrt{3}}{16\pi^2}\left(\frac{\sin^2\vartheta}{f_1^2}+\frac{16\cos^2\vartheta\sin^2\varphi}{25r_0^2(f_1')^2}\right).$$
Recalling that $\ds r_0^2(f_1')^2\leq \frac{1}{192\pi^2}$ due to (\ref{epsdel-ineq4bis}) and $\ds f_1^2\leq\frac{1}{\left(16\sqrt{3}(\pi^2+5)\right)^4}$ by (\ref{01-pinch})(\ref{epsdel-ineq3bis}), we obtain
$${\cal C}(m_0,V_0)(\xi,\nu)\geq\frac{1212\sqrt{3}}{16\pi^2}\left(\sin^2\vartheta+\cos^2\vartheta\sin^2\varphi\right)$$
which, combined with (\ref{auxil-trigo-bound2}), yields (\ref{strong-c-posit-ineq}) with $\varsigma=0.3$ at $(m_0,V_0,\xi,\nu)$.
\subsection{Concluding the proof of Proposition \ref{posit-conj-prop}}\label{conclu-conj-subsec}
By inspection of the smallness conditions (\ref{epsdel-ineq1}) through (\ref{epsdel-ineq8}) which $\varepsilon$ and $\delta_1$ must satisfy, we find that (\ref{epsdel-ineq7}) implies all others. Calculation yields the pinching:
$$1439\leq B_{120}+\frac{3}{2}B_{121}+2B_{110}B_{111}+\frac{36}{5\pi^2} \leq 1440$$
the right-hand side of which provides the condition:
$$\frac{\varepsilon}{2}+\delta_1\leq 2.96\times 10^{-15}$$
as a sufficient one for (\ref{epsdel-ineq7}), hence for all, to be satisfied. It leads us to take:
\begin{equation}\label{choice-eta1delta1}
\eta_1= 2.96\times 10^{-15},\ \delta_1= 1.48\times 10^{-15},
\end{equation}
in the statement of Proposition \ref{posit-conj-prop}. As for $\varsigma_1$, we choose the smallest value among the ones found along the way, namely: $\ds\varsigma_1=\frac{1}{4\pi^2}$.
\\
Finally, let us stress that the proof just completed obviously departs from that of \cite{firivi2} mentionned in Remark \ref{rem-firivi2}; in particular, in each of the above cases, the origin of the blow up rate (quadratic or cubic) chosen for the positive lower bound on the $c$-curvature can readily be traced back to the expression of ${\cal C}(m_0,V_0)(\xi,\nu)$ itself, relying on Lemma \ref{posit-minor-lem1} and Lemma \ref{mainperturb-lem}.
\section{$c$-curvature almost-positivity near the origin}\label{sec-zero}
In this section, we prove Theorem \ref{strong-c-posit-th} at $(m_0,V_0,\xi,\nu)$ when $\ds d\left( m_0,\exp_{m_0}(V_0)\right)$ is small.
\begin{Proposition}\label{posit-zero-prop}
There exists a triple of small (strictly) positive real numbers $(\eta_2,\delta_2,\varsigma_2)$ such that ${\cal C}(m_0,V_0)(\xi,\nu)$ satisfies the lower bound (\ref{strong-c-posit-ineq}) with $\varsigma=\varsigma_2$, provided $\varepsilon=|K-1|_{C^2(S)}\leq\eta_2$ and $|V_0|\leq \delta_2$.
\end{Proposition}
\textbf{Proof.} As already observed, we may take $V_0\not=0$ with no loss of generality. Dropping the subscript of $\delta_2$, we take $\ds\bar{r}_0\leq\frac{\pi}{2}$ by assuming $\varepsilon$ and $\delta$ (\ref{epsdel-ineq9})-small. We use the Maclaurin type approximation of $\ds \frac{f_1^3}{\bar{f}_1^3}\ {\cal C}(m_0,V_0)(\xi,\nu)$ obtained in Corollary \ref{c-curv-perturb-cor} and proceed to specify it further as $r_0\downarrow 0$. As regards its first summand, namely $\overline{{\cal C}}(m_0,V_0)(\xi,\nu)$, the expression (\ref{c-curv-sphere-eq}) prompts us to define constants $c_{11},\ldots,c_{14}$ as done in Appendix \ref{app-b}. These definitions imply at once that the absolute value of:
$$\overline{{\cal C}}(m_0,V_0)(\xi,\nu)-\frac{2\kappa\bar{r}_0^2}{45}\sin^2\vartheta\sin^2\varphi-\frac{2\kappa}{3}\left(1+\frac{2\bar{r}_0^2}{5}\right)\sin^2\vartheta\cos^2\varphi$$
$$-\frac{2\kappa}{3}\left(1+\frac{2\bar{r}_0^2}{15}\right)\cos^2\vartheta\sin^2\varphi-\frac{4\kappa}{3}\left(1+\frac{\bar{r}_0^2}{5}\right)\cos\vartheta\sin\vartheta\cos\varphi\sin\varphi$$
is bounded above by:
$$\kappa\bar{r}_0^3( c_{11}\sin^2\vartheta\sin^2\varphi+c_{12}\sin^2\vartheta\cos^2\varphi$$
$$+c_{13}\cos^2\vartheta\sin^2\varphi+c_{14}|\cos\vartheta\sin\vartheta\cos\varphi\sin\varphi|).$$
Let us now focus on the second summand, namely on the expression
$$E_4:=\frac{r_0\psi_2\sin^2\vartheta}{\bar{f_1}}\left( {\cal S}_{\bar{r}_0}(t)(1)-\frac{\bar{f}_0{\cal S}_{\bar{r}_0}(t^2)(1)}{\bar{f_1}}\right)$$
$$-\frac{2r_0\psi_0{\cal S}_{\bar{r}_0}(t^2-t)(1)}{\bar{f_1}}\left(\cos^2\vartheta - \cos^2(\vartheta+\varphi)\right)$$
$$-\frac{4r_0\psi_1\cos\vartheta\sin\vartheta\sin\varphi}{\bar{f_1}}\left( {\cal S}_{\bar{r}_0}(t)(1)-\frac{\bar{f}_0{\cal S}_{\bar{r}_0}(t^2)(1)}{\bar{f_1}}\right)$$
and rewrite, on the one hand:
$$\frac{r_0}{\bar{f_1}}\left( {\cal S}_{\bar{r}_0}(t)(1)-\frac{\bar{f}_0{\cal S}_{\bar{r}_0}(t^2)(1)}{\bar{f_1}}\right)$$
as: $\ds r_0 {\cal S}_{\bar{r}_0}(t-t^2)(1)+\sqrt{\kappa}r_0^2\left[A_1(\bar{r}_0) {\cal S}_{\bar{r}_0}(t)(1)-A_2(\bar{r}_0) {\cal S}_{\bar{r}_0}(t^2)(1)\right],$ where\footnote{so that: $A_1(\bar{r}_0)=\frac{1}{\bar{r}_0}\left(\frac{1}{\bar{f}_1}-1\right), A_2(\bar{r}_0)=\frac{1}{\bar{r}_0}\left(\frac{\bar{f}_0}{\bar{f}^2_1}-1\right)$}:
$$A_1(\tau):=\frac{\tau-\sin\tau}{\tau\sin\tau},\ A_2(\tau):=\frac{\tau^2\cos\tau-\sin^2\tau}{\tau\sin^2\tau}$$
(and note that two additional constants $c_{15},c_{16}$ are defined accordingly as in Appendix \ref{app-b}), on the other hand:
$$\frac{r_0}{\bar{f_1}} {\cal S}_{\bar{r}_0}(t^2-t)(1)=r_0{\cal S}_{\bar{r}_0}(t^2-t)(1)+\sqrt{\kappa}r_0^2 A_1(\bar{r}_0) {\cal S}_{\bar{r}_0}(t^2-t)(1).$$
Furthermore, the Maclaurin expansion of ${\cal S}_{\bar{r}_0}(t^2-t)(1)$ prompts us to write:
$$r_0{\cal S}_{\bar{r}_0}(t^2-t)(1)=-\frac{r_0}{12}+\kappa r_0^3 A_3(\bar{r}_0)$$
(defining so the auxiliary function $A_3$ and, accordingly, a constant $c_{17}$ as in Appendix \ref{app-b}). Gathering terms of same order and recalling (\ref{auxil-trigo-bound}), we obtain that the absolute value of:
{\small $$E_4-\frac{r_0}{6}\left[2\sin\vartheta\sin\varphi\sin(\vartheta-\varphi)\partial_1K(0)+\left(2\sin\vartheta\cos\varphi\sin(\vartheta-\varphi)+\sin^2(\vartheta-\varphi)\right)\partial_2K(0)\right]$$}
is bounded above by:
$$2\sqrt{\kappa}\varepsilon r_0^2 [\left( 8(c_{15}c_{6}+c_{16}c_{7}) +(c_{6}+c_{7})c_{15}\right)\sin^2\vartheta$$
$$+\left( 4(c_{15}c_{6}+c_{16}c_{7}) + 2 (c_{6}+c_{7})c_{15}\right)\cos^2\vartheta\sin^2\varphi ]$$
$$+2c_{17}\ \kappa\varepsilon r_0^3 \left(9\sin^2\vartheta+6\cos^2\vartheta\sin^2\varphi\right)\ .$$
Combining the latter inequality with the one derived above for the first summand $\overline{{\cal C}}(m_0,V_0)(\xi,\nu)$ of the expansion of ${\cal C}(m_0,V_0)(\xi,\nu)$ given in Corollary \ref{c-curv-perturb-cor}, we infer that, if we consider the decomposition:
$$\frac{f_1^3}{\bar{f}_1^3}\ {\cal C}(m_0,V_0)(\xi,\nu)=I+II+III$$
with
$$I:=\frac{\kappa}{3}\left(1+\frac{23\bar{r}_0^2}{30}\right)\sin^2\vartheta\cos^2\varphi+\frac{\kappa}{3}\left(1+\frac{\bar{r}_0^2}{10}\right)\cos^2\vartheta\sin^2\varphi$$
$$-\frac{2\kappa}{3}\left(1+\frac{2\bar{r}_0^2}{5}\right)\sin\vartheta\cos\vartheta\sin\varphi\cos\varphi\ ,$$
and
$$II:=\frac{\kappa}{3}\sin^2(\vartheta-\varphi)+\frac{\kappa\bar{r}_0^2 }{180}\left ( \sin^2\vartheta\cos^2\varphi+\cos^2\vartheta\sin^2\varphi+4\sin^2\vartheta\sin^2\varphi\right)$$
$$+\frac{r_0}{6}\ [2\sin\vartheta\sin\varphi\sin(\vartheta-\varphi)\partial_1K(0)$$
$$+\left(2\sin\vartheta\cos\varphi\sin(\vartheta-\varphi)+\sin^2(\vartheta-\varphi)\right)\partial_2K(0)]$$
and $III$ so defined, then the quantity:
$$\left| III-\frac{\kappa\bar{r}_0^2 }{180}\left( \sin^2\vartheta\cos^2\varphi+\cos^2\vartheta\sin^2\varphi+4\sin^2\vartheta\sin^2\varphi\right)\right|$$
is altogether bounded above by:
$$\varepsilon r_0^2\sin^2\vartheta\left( \frac{338C_1^3\pi^8}{\bar{f}_1^3}+2\sqrt{\kappa}\left[ 8(c_{15}c_{6}+c_{16}c_{7}) +(c_{6}+c_{7})c_{15}\right]\right)$$
$$+\varepsilon r_0^2\cos^2\vartheta\sin^2\varphi\left(\frac{268C_1^3\pi^8}{\bar{f}_1^3}+ 2\sqrt{\kappa}\left[ 4(c_{15}c_{6}+c_{16}c_{7}) + 2 (c_{6}+c_{7})c_{15}\right]\right)$$
$$+\varepsilon r_0^3 \ 2\kappa c_{17} \left(9\sin^2\vartheta+6\cos^2\vartheta\sin^2\varphi\right)$$
$$+\kappa\bar{r}_0^3\ ( c_{11}\sin^2\vartheta\sin^2\varphi+c_{12}\sin^2\vartheta\cos^2\varphi)$$
$$+\kappa\bar{r}_0^3\ (c_{13}\cos^2\vartheta\sin^2\varphi+c_{14}|\cos\vartheta\sin\vartheta\cos\varphi\sin\varphi|).$$
$$$$
Now, let us discuss separately the positivity of each summand $I,II,III$. Noting that
$$I\geq \frac{2\kappa}{3}|\cos\vartheta\sin\vartheta\cos\varphi\sin\varphi|\left( \sqrt{\left(1+\frac{\bar{r}_0^2}{10} \right)\left(1+\frac{23\bar{r}_0^2}{30} \right)}-\left(1+\frac{2\bar{r}_0^2}{5} \right)\right),$$
we find $I\geq 0$ provided $\ds \bar{r}_0\leq\frac{2}{\sqrt{5}}$ which holds if $\varepsilon$ and $\delta$ are (\ref{epsdel-ineq10})-small. Next, we have:
$$II\geq \frac{\kappa}{3}\sin^2(\vartheta-\varphi)+\frac{\kappa\bar{r}_0^2 }{180}\left ( \sin^2\vartheta+\sin^2\varphi+2\sin^2\vartheta\sin^2\varphi\right)$$
$$-\frac{\varepsilon r_0}{6}\ \left(4|\sin\vartheta|+|\sin(\vartheta-\varphi)|\right)|\sin(\vartheta-\varphi)|\ ,$$
hence
$$II \geq \frac{\kappa}{9}\sin^2(\vartheta-\varphi)+\frac{\kappa\bar{r}_0^2 }{360}\left ( \sin^2\vartheta+\sin^2\varphi\right)+\sin^2(\vartheta-\varphi)\left(\frac{\kappa}{9}-\frac{\varepsilon r_0}{6}\right)$$
$$+\frac{\kappa}{9}\sin^2(\vartheta-\varphi)+\frac{\kappa\bar{r}_0^2 }{360} \sin^2\vartheta-\frac{2\varepsilon r_0}{3}|\sin\vartheta\sin(\vartheta-\varphi)|\ .$$
So, assuming provisionally $III\geq 0$, and under the further smallness conditions\footnote{implied, for instance, by (\ref{epsdel-ineq3}) and (\ref{epsdel-ineq10})}:
$$\varepsilon\delta\leq\frac{2}{3},\ \varepsilon\leq \frac{1}{6\sqrt{10}},$$
the first of which implies $\ds \left(\frac{\kappa}{9}-\frac{\varepsilon r_0}{6}\right)\geq 0$, the second of which ensures that the second line of our last lower bound on $II$ is identically non-negative, we obtain:
\begin{equation}\label{c-curv-min-zero}
\frac{f_1^3}{\bar{f}_1^3}\ {\cal C}(m_0,V_0)(\xi,\nu)\geq \frac{\kappa}{9}\sin^2(\vartheta-\varphi)+\frac{\kappa^2}{360}\ r_0^2 \left( \sin^2\vartheta+\sin^2\varphi\right).
\end{equation}
From $r_0\leq\bar{r}_0\leq\frac{\pi}{2}$ combined with Remark \ref{sturm-rk}, we find $\ds\frac{f_1}{\bar{f}_1}\leq\sqrt{\kappa}$, with $\sqrt{\kappa}\leq 1+\frac{\varepsilon}{2}\leq1+\frac{1}{12\sqrt{10}}$ due to our last smallness assumption on $\varepsilon$. It yields $\ds\frac{f_1^3}{\bar{f}_1^3}\leq 1.1$ and the latter, plugged into (\ref{c-curv-min-zero}) proves Proposition \ref{posit-zero-prop} with $\ds \varsigma=\frac{1}{396}$ in (\ref{strong-c-posit-ineq}).
\medskip\\
Finally, let us discuss the non-negativity of $III$. From $\bar{r}_0\leq\frac{\pi}{2}$, we have $\ds \bar{f}_1(\bar{r}_0)\geq \frac{2}{\pi}$; moreover, we just saw that $\sqrt{\kappa}$ is bounded above by $\ds 1+\frac{1}{12\sqrt{10}}<1.027$. So the constants $C_{2},C_{3}$ defined in Appendix \ref{app-b} can be used as upper bounds on the coefficients respectively of $\varepsilon r_0^2\sin^2\vartheta$ and $\varepsilon r_0^2\cos^2\vartheta \sin^2\varphi$ in the lengthy expression which controls $\ds\left| III - \frac{\kappa\bar{r}_0^2 }{180}\ldots \right|$ (\textit{cf. supra}). Using them and recalling (\ref{auxil-trigo-bound}), we infer from the control just mentionned that:
\begin{eqnarray}
\frac{1}{r_0^2}III &\geq& \sin^2\vartheta\left[ \frac{1}{180}- \varepsilon( C_{2}+ 19c_{17}\delta) - \frac{115}{100}\delta (c_{11} +c_{12} +\frac{1}{2} c_{14})\right]\nonumber\\
&+& \cos^2\vartheta \sin^2\varphi\left[ \frac{1}{180}- \varepsilon( C_{3}+ 13c_{17}\delta) - \frac{115}{100}\delta (c_{13} +\frac{1}{2} c_{14})\right].\nonumber
\end{eqnarray}
Therefore $III\geq 0$ provided $\varepsilon$ and $\delta$ are taken (\ref{epsdel-ineq11})(\ref{epsdel-ineq12})-small. Proposition \ref{posit-zero-prop} is proved.
\paragraph{Concluding the proof of Proposition \ref{posit-zero-prop}.} By inspection of the smallness conditions (\ref{epsdel-ineq9}) through (\ref{epsdel-ineq12}) which $\varepsilon$ and $\delta_2$ must satisfy, we find that (\ref{epsdel-ineq11}) is the strongest one bearing on $\varepsilon$, because $C_2$ (like $C_3<C_2$) is $O(10^{18})$ while the constants $c_i$'s (with $11\leq i\leq 17$) are $O(1)$. It is also the strongest smallness condition on $\delta=\delta_2$ since setting $\varepsilon=1$ in (\ref{epsdel-ineq11}) yields $\ds\delta\leq\frac{1}{78}$. We will thus take:
\begin{equation}\label{choice-delta2}
\delta_2=0.01
\end{equation}
and, plugging this choice in (\ref{epsdel-ineq11}), get: $\ds \varepsilon C_2\leq 1.214\times 10^{-3}$. Since $C_2\leq 1.4\times 10^{18}$, it leads us to take:
\begin{equation}\label{choice-eta2}
\eta_2=8.6\times 10^{-22}.
\end{equation}
So, Proposition \ref{posit-zero-prop} holds with $(\eta_2,\delta_2)$ as just chosen and $\ds \varsigma_2=\frac{1}{396}$ (as found above).
\section{$c$-curvature almost-positivity elsewhere}\label{sec-else}
In this section, we prove Theorem \ref{strong-c-posit-th} at $(m_0,V_0,\xi,\nu)$ when $\exp_{m_0}(V_0)$ stays away from $m_0$ and $m_0^*$ as specified\footnote{sticking to the notations of Propositions \ref{posit-conj-prop} and \ref{posit-zero-prop}} in the:
\begin{Proposition}\label{posit-else-prop}
There exists a couple of small (strictly) positive real numbers $(\eta_3,\varsigma_3)$ such that ${\cal C}(m_0,V_0)(\xi,\nu)$ satisfies the lower bound (\ref{strong-c-posit-ineq}) with $\varsigma=\varsigma_3$, provided $\varepsilon=|K-1|_{C^2(S)}\leq\eta_3$ and
$\ds\frac{1}{2}\delta_2\leq |V_0| \leq \left(1-\frac{1}{2}\delta_1\right)\ell_0$.
\end{Proposition}

\noindent\textbf{Proof.} The following pinching holds:
$$\frac{1}{2}\delta_2\sqrt{1-\varepsilon}\leq\bar{r}_0\leq\pi\left(1-\frac{1}{2}\delta_1\right)\sqrt{1+\varepsilon}\ .$$
Recalling (\ref{choice-eta1delta1}) and assuming that $\varepsilon\leq\eta_2$, it implies the other one:
\begin{equation}\label{rbar-pinch-else}
\frac{49}{100}\delta_2\leq\bar{r}_0\leq \left(1-\frac{1}{4}\delta_1\right)\pi,
\end{equation}
the right-hand side of which yields the estimate:
\begin{equation}\label{invf1bar-estim}
\frac{1}{\bar{f}_1}\leq\frac{\pi}{\sin\left(\frac{\pi}{4}\delta_1 \right)},
\end{equation}
recorded here for later use. From Corollary \ref{c-curv-perturb-cor} combined with (\ref{obv-psi-bound}), (\ref{auxil-trigo-bound}), $r_0\leq\pi$ and $\kappa\geq 1$, we may write: 
\begin{eqnarray}\label{c-curv-min-else1}
\frac{f^3_1}{\bar{f}^3_1}{\cal C}(m_0,V_0)(\xi,\nu) &\geq& \frac{1}{\kappa}\overline{{\cal C}}(m_0,V_0)(\xi,\nu) - \frac{\varepsilon}{\bar{f}^3_1}\sin^2\vartheta\left(338C_1^3\pi^{10}+20\pi(c_6+c_7)\right)\nonumber\\
&-& \frac{\varepsilon}{\bar{f}^3_1}\cos^2\vartheta\sin^2\varphi\left(268C_1^3\pi^{10}+20\pi(c_6+c_7)\right).
\end{eqnarray}
The inequality:
\begin{equation}\label{f1/f1bar-maj1}
\frac{f_1}{\bar{f}_1}\leq \frac{\sqrt{\kappa}\sin r_0}{\sin \bar{r}_0},
\end{equation}
obvious from Remark \ref{sturm-rk}, will be used below to deal with the left-hand side of (\ref{c-curv-min-else1}). As for the term $\ds\frac{1}{\kappa}\overline{{\cal C}}(m_0,V_0)(\xi,\nu)$ occurring in the right-hand side of (\ref{c-curv-min-else1}), recalling its expression (\ref{c-curv-sphere-eq}), we split it into two summands, namely, the square:
$$\overline{S}_1=2\left( \sin\vartheta\cos\varphi\sqrt{\frac{\bar{r}_0^2-\sin^2\bar{r}_0}{\bar{r}_0\sin^3\bar{r}_0}}- \cos\theta\sin\varphi\sqrt{\frac{\bar{r}_0^2-\sin^2\bar{r}_0}{\bar{r}_0^3\sin\bar{r}_0}}\right)^2,$$
and the remaining part, equal to:
{\small $$\overline{S}_2=\sin^2\vartheta\sin^2\varphi \frac{h_1(\bar{r}_0)}{\bar{r}_0^2\sin^2\bar{r}_0}+ 8\sin^2\vartheta\cos^2\varphi \frac{\cos \frac{\bar{r}_0}{2}h_2(\frac{\bar{r}_0}{2})}{\bar{r}_0\sin^3\bar{r}_0} +8 \cos^2\vartheta\sin^2\varphi \frac{\cos \frac{\bar{r}_0}{2}h_2(\frac{\bar{r}_0}{2})}{\bar{r}_0^3\sin\bar{r}_0},$$}
where
$$h_1(\tau)=\tau^2+\tau\sin\tau\cos\tau-2\sin^2\tau\ ,$$
$$h_2(\tau)=(\tau+\sin\tau\cos\tau)\sin\tau-2\tau^2\cos\tau\ .$$
Obviously, setting:
$$\mu_1(\tau):=\min\left(\frac{h_1(\tau)}{\tau^2\sin^2\tau}, \frac{8\cos \frac{\tau}{2}h_2(\frac{\tau}{2})}{\tau\sin^3\tau}, \frac{8\cos \frac{\tau}{2}h_2(\frac{\tau}{2})}{\tau^3\sin\tau}\right),$$
we have: $\ds\overline{S}_2\geq\mu_1(\bar{r}_0)(\sin^2\vartheta+ \cos^2\vartheta\sin^2\varphi )$, so we focus on a positive lower bound on $\mu_1(\bar{r}_0)$.\\
\indent To proceed further, let us distinguish two cases and split the proof accordingly.
\subsection{First case: $\ds\frac{49}{100}\delta_2\leq \bar{r}_0 \leq 1$}
In that case, on the one hand we may write:
\begin{equation}\label{invf1bar-maj-zero}
\frac{1}{\bar{f}_1}\leq\frac{1}{\sin 1}\ ,
\end{equation}
on the other hand, combining (\ref{f1/f1bar-maj1}) with the alternating series test applied to the Maclaurin series of $\sin\bar{r}_0$, we get:
$$\frac{f_1}{\bar{f}_1}\leq \frac{1}{1-\frac{\bar{r}_0^2}{6}}\leq \frac{6}{5}\ ,$$
so (\ref{c-curv-min-else1}) implies:
\begin{equation}\label{c-curv-min-else2}
\left(\frac{6}{5}\right)^3{\cal C}(m_0,V_0)(\xi,\nu) \geq \frac{1}{\kappa}\overline{{\cal C}}(m_0,V_0)(\xi,\nu) - \frac{\varepsilon C_4}{\sin^31}(\sin^2\vartheta+\cos^2\vartheta\sin^2\varphi)\ ,
\end{equation}
where $C_4$ is the constant defined in Appendix \ref{app-sec-else}. Besides, the following lemma holds:
\begin{Lemma}\label{htools-lem}
The function $h_1$ (resp. $h_2$) is increasing on $[0,\pi]$ (resp. on $\left[0,\frac{\pi}{2}\right]$). Furthermore, for each $\tau\in[0,1]$, the alternating series test holds for the Maclaurin series of $h_1(\tau)$ and $h_2(\tau)$, implying the lower bounds:
$$h_1(\tau)\geq\frac{2}{315}\tau^6(7-\tau^2),\
h_2(\tau)\geq\frac{4}{5}\tau^6\left(\frac{2}{9} - \frac{1}{21}\tau^2\right).$$
\end{Lemma}
The proof is lengthy but rather elementary hence left as an exercise. Combining this lemma with the standard bounds $\ds\sin\tau\leq\tau,\ \cos\frac{\tau}{2}\geq 1-\frac{\tau^2}{8}$, we get:
$$\mu_1(\bar{r}_0)\geq\min\left( \frac{14}{315}\bar{r}_0^2\left( 1- \frac{1}{7} \bar{r}_0^2\right),  \frac{1}{45}\bar{r}_0^2\left( 1- \frac{5}{28}\bar{r}_0^2\right)\right)\equiv  \frac{1}{45} \bar{r}_0^2\left( 1- \frac{5}{28}\bar{r}_0^2 \right)$$
and, from (\ref{rbar-pinch-else}), conclude: $\ds \mu_1(\bar{r}_0)\geq 4.38\times 10^{-3}\ \delta_2^2$. This lower bound combined with (\ref{c-curv-min-else2}) implies:
$${\cal C}(m_0,V_0)(\xi,\nu) \geq \left(\frac{5}{6}\right)^3\ 2.19\times 10^{-3}\ \delta_2^2\ (\sin^2\vartheta+\cos^2\vartheta\sin^2\varphi)\ ,$$
provided $\varepsilon$ is taken (\ref{eps-ineq-else1})-small. Recalling (\ref{auxil-trigo-bound2}), we thus obtain (\ref{strong-c-posit-ineq}) at $(m_0,V_0,\xi,\nu)$ with $\ds\varsigma=\frac{1}{4\pi^2}\left(\frac{5}{6}\right)^3\ 2.19\times 10^{-3}\ \delta_2^2$. Here, the value of $\delta_2$ is the one chosen at the end of Section \ref{sec-zero}, namely $\delta_2=0.01$; plugging it in the preceding formula, and in (\ref{eps-ineq-else1}) together with the sharp bound $C_4\leq3.6\times 10^{18}$, leads us to take:
\begin{equation}\label{eta3sigma3choice1}
\eta_3\leq 3.62\times 10^{-26},\ \varsigma_3= 3.21\times 10^{-9}.
\end{equation}

\subsection{Second case: $\ds 1\leq \bar{r}_0 \leq \left(1-\frac{\delta_1}{4}\right)\pi$}
Back to (\ref{c-curv-min-else1}), using (\ref{invf1bar-estim}), we now have:
$$\frac{f^3_1}{\bar{f}^3_1}{\cal C}(m_0,V_0)(\xi,\nu) \geq \frac{1}{\kappa}\overline{{\cal C}}(m_0,V_0)(\xi,\nu) - \frac{\varepsilon\pi^3C_{4}}{\sin^3\left(\frac{\delta_1}{4}\pi\right)}(\sin^2\vartheta+\cos^2\vartheta\sin^2\varphi)\ ,$$
hence:
\begin{equation}\label{c-curv-min-else3}
\frac{f^3_1}{\bar{f}^3_1}{\cal C}(m_0,V_0)(\xi,\nu) \geq \left( \mu_1(\bar{r}_0) - \frac{\varepsilon\pi^3C_{4}}{\sin^3\left(\frac{\delta_1}{4}\pi\right)}\right) (\sin^2\vartheta+\cos^2\vartheta\sin^2\varphi)\ .
\end{equation}
Moreover, using:
$$\sin\tau=\sin(\pi-\tau)\leq\pi-\tau,\ {\rm and}\ \cos\frac{\tau}{2}=\sin\frac{(\pi-\tau)}{2}\geq\frac{(\pi-\tau)}{2}-\frac{(\pi-\tau)^3}{48},$$
with $\tau=\bar{r}_0$, we have:
$$\frac{\cos\frac{\bar{r}_0}{2}}{\sin\bar{r}_0}\geq\frac{1}{2}-\frac{(\pi-1)^2}{48}\geq 0.4\ ,$$
therefore:
$$\frac{8\cos \frac{\bar{r}_0}{2}h_2(\frac{\bar{r}_0}{2})}{\bar{r}_0^3\sin\bar{r}_0}\geq\frac{3.2}{\pi^3}\ h_2\left(\frac{1}{2}\right)\geq 2.6\times 10^{-4}.$$
Besides, we directly get:
$$\frac{h_1(\bar{r}_0)}{\bar{r}_0^2\sin^2\bar{r}_0}\geq\frac{h_1(1)}{\pi^2}\geq 3.9\times 10^{-3},$$
and thus conclude: $\ds \mu_1(\bar{r}_0)\geq 2.6\times 10^{-4}$. Finally, from (\ref{f1/f1bar-maj1}), we infer the bound:
$$\frac{f_1}{\bar{f}_1}\leq \left(1+\frac{\varepsilon}{2}\right)\frac{\sin r_0}{\sin \bar{r}_0}\ ,$$
and from the identity:
$$\sin r_0\equiv\sin\bar{r}_0\ \cos[(\sqrt{\kappa}-1)r_0]-\cos\bar{r}_0\ \sin[(\sqrt{\kappa}-1)r_0]\ ,$$
we readily get:
$$\frac{\sin r_0}{\sin \bar{r}_0}\leq 1+\frac{\sin[(\sqrt{\kappa}-1)r_0]}{\sin\left(\frac{\delta_1}{4}\pi \right)}\leq 1+\frac{2\varepsilon}{\delta_1\left(1-\frac{\delta_1^2\pi^2}{96}\right)}\ .$$
It prompts us to take $\varepsilon$ (\ref{eps-ineq-else2})-small with $\delta_1$ as chosen at the end of Section \ref{sec-conj} (namely $\delta_1=1.48\times 10^{-15}$), in order to keep the ratio $\ds\frac{f_1}{\bar{f}_1}$ below $\frac{6}{5}$.
Plugging in (\ref{c-curv-min-else3}) the latter upper bound together with the former lower bound on $\mu_1(\bar{r}_0)$, we obtain:
$${\cal C}(m_0,V_0)(\xi,\nu) \geq \left(\frac{5}{6}\right)^3 1.3\times 10^{-4}\ (\sin^2\vartheta+ \cos^2\vartheta\sin^2\varphi )\ ,$$
provided $\varepsilon$ is taken (\ref{eps-ineq-else3})-small. Recalling (\ref{auxil-trigo-bound2}), we conclude that ${\cal C}(m_0,V_0)(\xi,\nu)$ satisfies (\ref{strong-c-posit-ineq}) with:$$\varsigma\leq \frac{1}{4\pi^2}\left(\frac{5}{6}\right)^3 1.3\times 10^{-4},$$
so here, it is sufficient to take: $\ds\varsigma_3\leq 1.9\times 10^{-6}$, a condition well satisfied by the value chosen in (\ref{eta3sigma3choice1}) for $\varsigma_3$. Finally, recalling that $\delta_1$ was taken equal to $1.48\times 10^{-15}$, the smallness condition (\ref{eps-ineq-else3}) on $\varepsilon$ leads us to take:
\begin{equation}\label{eta3choice2}
\eta_3= 1.8\times 10^{-69}.
\end{equation}
This tiny value (compare with (\ref{choice-eta1delta1})(\ref{choice-eta2})(\ref{eta3sigma3choice1})) reflects the fact that a perturbation device from the constant curvature case becomes outrageously rough as $|V_0|\uparrow \ell_0$ (\textit{i.e.} getting close to the first conjugate point).
\section{Proof of Theorem \ref{strong-c-posit-th}}\label{sec-synth}
The proof of Theorem \ref{strong-c-posit-th} at $(m_0,V_0,\xi,\nu)$ goes by combining Propositions \ref{posit-conj-prop}, \ref{posit-zero-prop} and \ref{posit-else-prop}. Doing so, we first observe that the assumption made on $|V_0|$ in Proposition \ref{posit-else-prop} overlaps, as it should, the corresponding ones of Propositions \ref{posit-conj-prop} and \ref{posit-zero-prop}. Next, since $\varepsilon$ should now fulfill \emph{all} the smallness conditions stated on it in Sections \ref{sec-conj}, \ref{sec-zero} and \ref{sec-else}, we take $\eta$ in the statement of Theorem \ref{strong-c-posit-th} equal to:
$$\eta=\min(\eta_1,\eta_2,\eta_3)\equiv \eta_3= 1.8\times 10^{-69}.$$
Similarly, we choose:
$$\varsigma=\min(\varsigma_1,\varsigma_2,\varsigma_3)\equiv \varsigma_3 = 3.21\times 10^{-9}.$$
\appendix
\section{Proof of Lemma \ref{mainperturb-lem}}\label{app-a}
We will proceed stepwise in the Fermi chart along $V_0$, using repeatedly the Maclaurin theorem, the
solution map
${\cal S}_{\bar{r}_0}$ and its contraction property, to derive estimates
at $((0,r_0),t)$, uniform in $t\in[0,1]$, on the expressions appearing in the conclusion of Lemma \ref{mainperturb-lem} and also on
$³|D^2X|$ and
$|D^j{\cal K}|$ for
$j=1,2$, where  ${\cal K}=K\circ X$.
\subsection{Estimates of order 0}  
\subsubsection{Basic estimates}
From Remark \ref{sturm-rk}, we may take $B_{101}=1$. 
Besides, we have:
\begin{equation}\label{0-perturK-estim}
\Vert \kappa -{\cal K} \Vert \leq\varepsilon\ \min(1, r_0)\ .
\end{equation}
On the axis of the Fermi chart, the functions $\tilde{f}_a=f_a-\bar{f}_a$ (with $a=0,1$) satisfy:
$$\frac{d^2\tilde{f}_a}{dt^2}+r_0^2\kappa\ \tilde{f}_a=\phi_{0a}\ {\rm with}\ \phi_{0a}=
r_0^2(\kappa -{\cal K})f_a\ .$$
Combining the latter with
(\ref{contraction-estim}) applied to ${\cal S}_{\bar{r}_0}$, and (\ref{0-perturK-estim}), we get:
$$\Vert\tilde{f}_a\Vert\leq\frac{\varepsilon}{2}\min( r_0^2, r_0^3)\ \Vert f_a\Vert .$$
If $a=0$, since $\Vert f_0\Vert\leq\Vert \tilde{f_0}\Vert + \Vert \bar{f_0}\Vert \leq \Vert \tilde{f_0}\Vert + 1$, we infer:
$$\Vert \tilde{f_0}\Vert\leq\frac{\mu}{1-\mu}\ {\rm with}\ \mu=\frac{\varepsilon}{2}\min( r_0^2, r_0^3),$$
while if $a=1$, recalling Remark \ref{sturm-rk}, we get at once: $\Vert \tilde{f_1}\Vert\leq \mu$. Since $\ds\varepsilon\leq \frac{1}{\pi^2}$,
we have $\varepsilon r_0^2\leq 1$ (recalling Remark \ref{comparison-rk}), an inequality used throughout this appendix. So we readily obtain:
$$\Vert \tilde{f_0}\Vert\leq\varepsilon\ \min\left( r_0^2, \frac{ r_0^3}{2-\varepsilon r_0^3} \right) .$$
In particular, regarding the first line of conclusion of the lemma for $k=0$, we may take $\ds B_{200}=1,B_{201}=\frac{1}{2}$, which yields
$B_{100}=2$ after use of the triangle
inequality. Similarly, setting $h_0=1$ and $h_1=t$, we find on the axis:
$$f_a-h_a={\cal S}_{\bar{r_0}}\left(-r_0^2\ {\cal K}\ h_a+r_0^2(\kappa -{\cal K})(f_a-h_a)\right)$$
 for $a=0,1$. Combining (\ref{contraction-estim}) with an argument as the one above for $\tilde{f}_0$ yields:
$$\Vert f_a-h_a\Vert\leq \frac{r_0^2\Vert {\cal K}\Vert}{2-\varepsilon r_0^2}\leq r_0^2\Vert {\cal K}\Vert$$
hence the inequalities:
\begin{equation}\label{0-basineq}
\Vert f_0-1\Vert\leq r_0^2(1+\varepsilon),\ \Vert f_1-t\Vert\leq r_0^2(1+\varepsilon)\ ,
\end{equation}
recorded here for later use.
\subsubsection{Estimates on Maclaurin approximations}
The first order Maclaurin approximation of ${\cal K}$ at $t=0$ satisfies the estimate:
\begin{equation}\label{0-mclaurK-estim}
\Vert{\cal K}-\kappa-t r_0\ \partial_2K(0)\Vert\leq \frac{1}{2}\varepsilon\ r_0^2\ .
\end{equation}
The latter combined with the triangle inequality is used to evaluate the remainder of the first non trivial Maclaurin approximation
of $\phi_{0a}$ at $t=0$,
namely of
$\phi_{0a}+t^{a+1}r_0^3\ \partial_2K(0)$
written as:
$$\phi_{0a}+t^{a+1}r_0^3\ \partial_2K(0)=-r_0^2\left({\cal K}-\kappa-t r_0\ \partial_2K(0)\right) f_a+t r_0^3\ \partial_2K(0)\ (h_a-f_a).$$
It leads us to the upper bound:
$$\Vert \phi_{0a}+t^{a+1}r_0^3\ \partial_2K(0)\Vert\leq \frac{1}{2}\varepsilon\ r_0^4\ \Vert f_a\Vert+\varepsilon\ r_0^3\ \Vert f_a-h_a\Vert$$
which, combined with (\ref{0-basineq}) and (\ref{contraction-estim}), yields
for
$$\tilde{f}_a+r_0^3\ \psi_0\ {\cal S}_{\bar{r_0}}(t^{a+1})\equiv{\cal S}_{\bar{r_0}}\left(\phi_{0a}+t^{a+1}r_0^3\ \partial_2K(0)\right)$$
the desired second line of conclusion with $\ds B_{30a}=\frac{1}{4}B_{10a}+\frac{\pi}{2}\left(1+\frac{1}{\pi^2}\right)$.
\subsection{Estimates of order 1}
\subsubsection{Basic estimates}
From the definition of ${\cal K}$ and $f_1$, we have at $(v_0,t)$:
$$D_1{\cal K}=f_1(t)\ (\partial_1K)(0,tr_0) ,\ D_2{\cal
K}=t\ (\partial_2K)(0,tr_0) .$$
Recalling Remark \ref{sturm-rk}, we conclude:
\begin{equation}\label{1-perturK-estim}
\forall i=1,2,\forall t\in[0,1], \ |D_i{\cal K}(v_0,t)|\leq\varepsilon,\ {\rm thus}\ \Vert D_{\nu}{\cal K}(v_0,.)\Vert\leq \sqrt{2}\varepsilon.
\end{equation}
Besides, if we apply $D_{\nu}$ to the Jacobi
equations:
\begin{equation}\label{fermijacob-eq}
\ddot{f}+|v|^2{\cal K}(v,t)f=0\ {\rm and}\ \ddot{f}+|v|^2\kappa f=0 ,
\end{equation}
then let $v=v_0=(0,r_0)$, we readily infer for $\tilde{f}_a$ the equation (still abbreviating freely $D_{\nu}$ by a prime):
$\ds\frac{d^2\tilde{f}'_a}{dt^2}+r_0^2\kappa\ \tilde{f}'_a=\phi_{1a}$, with:
$$\phi_{1a}=r_0^2(\kappa-{\cal K})f_a'-2r_0\cos\varphi\ \kappa \tilde{f}_a+2r_0\cos\varphi\ f_a(\kappa-{\cal
K})-r_0^2{\cal K}'f_a,$$
and for $\bar{f}_a$ the equation: $\ds\frac{d^2\bar{f}'_a}{dt^2}+r_0^2\kappa\ \bar{f}'_a=-2r_0\cos\varphi\ \kappa\bar{f}_a$. Recalling
(\ref{contraction-estim}), we get from the latter the auxiliary bound:
\begin{equation}\label{1-barf-auxil}
\Vert\bar{f}'_a\Vert\leq r_0\kappa\leq c_1
\end{equation}
(see Appendix \ref{app-b}), and from the former:
$$\Vert\tilde{f}'_a\Vert\leq \frac{1}{2}r_0^2\Vert\kappa-{\cal
K}\Vert\left(\Vert\tilde{f}'_a\Vert + \Vert\bar{f}'_a\Vert\right)+r_0\kappa\ \Vert \tilde{f}_a\Vert+r_0\Vert\kappa-{\cal K}\Vert\ \Vert
f_a\Vert+
\frac{1}{2}r_0^2\Vert {\cal K}'\Vert\ \Vert f_a\Vert,$$
after use of the triangle inequality. Previous bounds, namely (\ref{0-perturK-estim})(\ref{1-perturK-estim})(\ref{1-barf-auxil}) and those of
Lemma
\ref{mainperturb-lem} for
$k=0$, yield:
$$\Vert\tilde{f}'_a\Vert\leq\frac{1}{1-\frac{1}{2}\varepsilon r_0^2}\left(\frac{1}{2}\varepsilon\kappa r_0^3+B_{20a}\varepsilon\kappa
r_0^3+B_{10a}\varepsilon r_0+\frac{1}{2}B_{10a}\sqrt{2}\varepsilon r_0^2
\right)$$
hence the conclusion of the first line of the lemma holds for $k=1$ with:
$$B_{21a}=1+\pi^2+2B_{20a}(1+\pi^2)+B_{10a}(2+\pi\sqrt{2}),$$
and, combining the triangle inequality with the auxiliary bound on $\bar{f}'_a$, with:
$$B_{11a}=\pi +\frac{1}{\pi}(1+B_{21a}).$$
\subsubsection{Estimates on Maclaurin approximations}
From the expression found above for $D{\cal K}(v_0,t)$, we may write:
$$D_{\nu}{\cal K}(v_0,t)=t\partial_{\nu}K(0,tr_0)+\sin\varphi\ \partial_1K(0,tr_0)\ (f_1-t).$$
So, using the straightforward bound: $|\partial_{\nu}K(0,tr_0)-\partial_{\nu}K(0)|\leq\varepsilon r_0$ combined with the triangle inequality
and (\ref{0-basineq}), we obtain:
\begin{equation}\label{1-mclaurK-estim}
\Vert {\cal K}'-t\partial_{\nu}K(0)\Vert\leq (1+c_1)\varepsilon r_0\ .
\end{equation}
We wish now to estimate the remainder of the first non trivial Maclaurin approximation
of $\phi_{1a}$ at $t=0$,
namely the $\Vert .\Vert$ norm of the expression:
$$\phi_{1a}+2t\ h_a r_0^2\cos\varphi\ \partial_2K(0)+t\ h_ar_0^2\ \partial_{\nu}K(0).$$
To do so, we recast the latter as follows:
$$\ds=-2\kappa r_0\cos\varphi\ \tilde{f}_a-2r_0\cos\varphi\ f_a\left({\cal K}-\kappa-t r_0\ \partial_2K(0)\right)+2tr_0^2\cos\varphi\ 
\partial_2K(0)(h_a-f_a)$$
$$+r_0^2(\kappa -{\cal K})(\tilde{f}_a'+\bar{f}_a')-r_0^2\left({\cal K}'-t\partial_{\nu}K(0)\right)f_a-tr_0^2\
\partial_{\nu}K(0)(f_a-h_a)$$
and apply the triangle inequality combined with
(\ref{0-perturK-estim})(\ref{0-basineq})(\ref{0-mclaurK-estim})(\ref{1-barf-auxil})(\ref{1-mclaurK-estim}) and the bounds of the lemma on the
$\Vert .\Vert$ norms of $\tilde{f}_a,\tilde{f}_a'$. Observing that, if we apply the map ${\cal S}_{\bar{r_0}}$ to the preceding
expression and use (\ref{auxil-sol-formulas}), we recover
$\tilde{f}_a'+r_0^2\
\psi_1\ {\cal S}_{\bar{r_0}}(t^{a+1})$, and recalling (\ref{contraction-estim}), we infer that $\Vert\tilde{f}_a'+r_0^2\
\psi_1\ {\cal S}_{\bar{r_0}}(t^{a+1})\Vert$ is bounded above by:
$$\varepsilon r_0^3\left(
\frac{1}{2}B_{10a}\left(2+\pi+\frac{1}{\pi}\right)+B_{20a}\kappa+B_{21a}\frac{\varepsilon}{2}+
\frac{3}{2}r_0(1+\varepsilon)+\frac{\kappa}{2}\right).$$
The second line of conclusion of Lemma
\ref{mainperturb-lem} for $k=1$, indeed, follows with:
$$B_{31a}=\frac{1}{2}B_{10a}\left(2+\pi+\frac{1}{\pi}\right)+B_{20a}\left(1+\frac{1}{\pi^2}\right)+B_{21a}\frac{1}{2\pi^2}+
\frac{(3\pi+1)}{2}\left(1+\frac{1}{\pi^2}\right).$$
\subsection{Estimates of order 2}
\subsubsection{Basic estimates}
As in \cite{delge}, applying twice $D_{\nu}$ to the geodesic equation with initial conditions $(0,v)$, then letting $v=v_0=(0,r_0)$, and
recalling the 2-dimensional formulas given after Definition \ref{fermi-chart-def} for the derivatives of the Christoffel symbols on the axis of
the Fermi chart, yields for
$D_{\nu\nu}X^i(t)=D_{\nu\nu}X^i(0,v_0,t)$ the following equations, with zero initial conditions:
\begin{eqnarray}
\frac{d^2}{dt^2}\left(D_{\nu\nu}X^1\right)+r_0^2{\cal K}\ D_{\nu\nu}X^1= &-& 4r_0\cos\varphi\sin\varphi\ {\cal K}f_1\nonumber\\ &-& r_0^2\sin^2\varphi\
f_1^2\ \left((\partial_1K)\circ X\right)\nonumber\\
&-& 2r_0^2\sin\varphi\cos\varphi\  tf_1\left((\partial_2K)\circ X\right),\nonumber
\end{eqnarray}
$$\frac{d^2}{dt^2}\left(D_{\nu\nu}X^2\right)=4r_0\sin^2\varphi\ {\cal K} f_1\dot{f}_1+ r_0^2\sin^2\varphi\ f_1^2\ \left((\partial_2K)\circ
X\right).$$
To treat the first equation, we view ${\cal K}$ as a perturbation of $\kappa$ and apply the solution map ${\cal
S}_{\bar{r}_0}$ and the estimates (\ref{contraction-estim}) (\ref{0-perturK-estim}) and that on $\Vert f_1\Vert$; to treat the second equation,
we use our estimates on $\Vert{\cal K}\Vert$ and $\Vert f_1\Vert$ and note the further one:
$$\dot{f}_1(t)=1+\int_0^t\ddot{f}_1(\theta)d\theta\equiv 1-r_0^2\int_0^t{\cal
K}(\theta)f_1(\theta) d\theta\Longrightarrow\Vert\dot{f}_1\Vert\leq1+r_0^2(1+\varepsilon)\leq 2+\pi^2.$$ We readily find:
\begin{equation}\label{D2X-estim}
\Vert D_{\nu\nu}X^1\Vert\leq c_2r_0,\  \Vert D_{\nu\nu}X^2\Vert\leq c_3r_0.
\end{equation}\\
Next, we calculate the expression of $D_{\nu\nu}{\cal K}(v_0,t)$ and obtain:
\begin{eqnarray}
D_{\nu\nu}{\cal K}(v_0,t) &=& \partial_1K(0,tr_0)\ D_{\nu\nu}X^1+\partial_2K(0,tr_0)\ D_{\nu\nu}X^2\nonumber\\
& & +\partial_{11}K(0,tr_0)\ f_1^2\sin^2\varphi+ 2\partial_{12}K(0,tr_0)\ tf_1\sin\varphi \cos\varphi\nonumber\\
& & +\partial_{22}K(0,tr_0)\ t^2\cos^2\varphi\ ,\nonumber
\end{eqnarray}
from what we infer, using (\ref{D2X-estim}) combined with Remark \ref{sturm-rk}:
\begin{equation}\label{2-perturK-estim}
\Vert D_{\nu\nu}{\cal K}(v_0,.)\Vert \leq (c_2+c_3)\varepsilon r_0+2\varepsilon\leq c_4\varepsilon.
\end{equation}\\
Now, we apply $D_{\nu\nu}$ to (\ref{fermijacob-eq}) and get, on the one hand:
$$\frac{d^2}{dt^2}\left(\bar{f}_a''\right)+r_0^2\kappa \bar{f}_a''=-2\kappa \bar{f}_a-4\kappa r_0\cos\varphi\ \bar{f}_a'\ ,$$
from what, recalling (\ref{1-barf-auxil}), we infer the auxiliary bound:
\begin{equation}\label{2-barf-auxil}
\Vert\bar{f}_a''³\Vert\leq c_5\ ,
\end{equation}
on the other hand:
$$\frac{d^2}{dt^2}\left(\tilde{f}_a''\right)+r_0^2\kappa \tilde{f}_a''=\phi_{2a}\ ,$$
with:
\begin{eqnarray}
\phi_{2a} &=&r_0^2(\kappa-{\cal K})f_a''+2(\kappa\bar{f}_a-{\cal K}f_a)+4r_0\cos\varphi(\kappa\bar{f}_a'-{\cal K}f_a')\nonumber\\
& & -4 r_0\cos\varphi\ f_aD_{\nu}{\cal K}-2r_0^2f_a'D_{\nu}{\cal K}-r_0^2f_aD_{\nu\nu}{\cal K}.\nonumber
\end{eqnarray}
Finally, from (\ref{contraction-estim}) applied (with $\omega=\bar{r}_0$) to the latter equation, we routinely derive the first line of
conclusion of Lemma \ref{mainperturb-lem} for $k=2$ with:
$$B_{12a}=c_5+\frac{1}{\pi^2}B_{22a}\ ,$$
after use of the triangle inequality combined with (\ref{2-barf-auxil}), and:
$$B_{22a}=6+\pi^2(4+c_5)+(4\sqrt{2}+\pi c_4)\pi B_{10a}+2\sqrt{2}\pi^2B_{11a}+2(1+\pi^2)\left( B_{20a}+2B_{21a} \right),$$
after use of the triangle inequality combined with (\ref{0-perturK-estim})(\ref{1-perturK-estim})(\ref{1-barf-auxil})(\ref{2-perturK-estim})
and the bounds of the same line of conclusion for $k=0,1$.
\subsubsection{Estimates on Maclaurin approximations}
Finally, in order to estimate the $\Vert .\Vert$ norm of $\tilde{f}_a''+r_0\psi_2{\cal S}_{\bar{r}_0}(t^{a+1})$, we note that the latter
is equal to ${\cal S}_{\bar{r}_0}\left(\phi_{2a}+2th_ar_0\partial_2K(0)+4r_0\cos\varphi\ th_a\partial_{\nu}K(0)\right)$, we recast the
argument of ${\cal S}_{\bar{r}_0}$ as follows:
$$\phi_{2a}+2th_ar_0\partial_2K(0)+4r_0\cos\varphi\ th_a\partial_{\nu}K(0) = r_0^2(\kappa-{\cal K})f_a''-2r_0^2f_a'D_{\nu}{\cal
K}-r_0^2f_aD_{\nu\nu}{\cal
K}$$
$$-2\kappa\tilde{f}_a-2tr_0 \partial_2K(0)(f_a-h_a) - 2\left( {\cal K}-\kappa-t
r_0\partial_2K(0)\right)f_a-4r_0\cos\varphi\ \kappa\tilde{f}_a'$$
$$+4r_0\cos\varphi(\kappa-{\cal K})f_a'-4r_0\cos\varphi\ f_a\left(D_{\nu}{\cal
K}-t\partial_{\nu}K(0)\right)+4tr_0\cos\varphi(h_a-f_a)\partial_{\nu}K(0)\ ,$$
and we apply (\ref{contraction-estim}) with $\omega=\bar{r}_0$ to the right-hand expression, combined with the triangle inequality, the
previous bounds of Lemma \ref{mainperturb-lem} and (\ref{0-perturK-estim})
(\ref{0-basineq})(\ref{0-mclaurK-estim})(\ref{1-perturK-estim})(\ref{1-mclaurK-estim})(\ref{2-perturK-estim}). Doing so term by term, we
obtain the second line of conclusion of Lemma \ref{mainperturb-lem} for $k=2$ with:
\begin{eqnarray}
B_{32a} &=& \frac{1}{2}B_{12a}+\sqrt{2}B_{11a}+\frac{1}{2}B_{10a}c_4+\left(1+\frac{1}{\pi^2}\right)B_{20a}+c_1\nonumber\\
& &+\frac{1}{2}B_{10a}+2\left(1+\frac{1}{\pi^2}\right)B_{21a}+2B_{11a}+2B_{10a}(1+c_1)+2c_1\nonumber\\
& \equiv & 3 c_1 + \frac{1}{2} (5+4c_1+c_4)B_{10a}+(\sqrt{2}+2)B_{11a}+\frac{1}{2}B_{12a}\nonumber\\
& &+\left(1+\frac{1}{\pi^2}\right)(B_{20a}+2B_{21a})\nonumber\ .
\end{eqnarray}
\section{Auxiliary universal constants and conditions}\label{app-b}
\subsection{List of constants for Section \ref{sec-perturb}}
$$\ds c_1=\pi+\frac{1}{\pi},\ c_2=4\left(\frac{3}{4\pi}+\frac{1}{\pi^2}+1\right),\ c_3=\frac{1}{2\pi}+2\left(1+\frac{1}{\pi^2}\right)\left(2+\pi^2\right),$$
$$c_4=\frac{11}{2}+10\pi+\frac{8}{\pi}+2\pi^3,\ c_5=1+\frac{1}{\pi^2}+2c_1^2,$$
$$c_{6}=\sup_{\tau\in\left[0,2\pi\right]}\left| \frac{\tau-\sin\tau}{\tau^3}\right|,\ c_{7}=\sup_{\tau\in\left[0,2\pi\right]}\left| \frac{\tau^2+2(\cos\tau-1)}{\tau^4}\right|,$$
$$c_{8}=\sup_{\tau\in\left[0,2\pi\right]}\left| \frac{\tau\cos\tau-\sin\tau}{\tau^2}\right|,\ c_{9}=\sup_{\tau\in\left[0,2\pi\right]}\left| \frac{\tau\cos\tau-\sin\tau}{\tau^3}\right|,$$
$$c_{10}=\sup_{\tau\in\left[0,2\pi\right]}\left| \frac{\cos\tau\sin\tau-\tau}{\tau^3}\right|,$$
$$C_1=\max\left(\max_{a=0,1;k=0,1,2}\left(B_{1ka},B_{3ka}\right),8c_6,8c_7,\frac{19}{18}c_8,\frac{10}{9}c_9,c_{10} \right).$$
\subsection{List of conditions on $\varepsilon$ and $\delta=\delta_1$ for Section \ref{sec-conj}}\label{app-sec-conj}
\begin{equation}\label{epsdel-ineq1}
\frac{|\cos\varphi|}{\pi}\left(1-\frac{\frac{\varepsilon}{2}}{1-\left(\frac{\varepsilon}{2}+\delta\right)}-\frac{\pi^2}{2}\left(\frac{\varepsilon}{2}+\delta\right)^2 \right)-\varepsilon\pi B_{211}>0
\end{equation}
an inequality to be used only in subsection \ref{subsec-cosnot0} with $|\cos\varphi|$ replaced by $\ds\frac{1}{7704}$;
\begin{equation}\label{epsdel-ineq2}
\frac{2\varepsilon}{1-\left(\frac{\varepsilon}{2}+\delta\right)}+2\pi^2\left(\frac{\varepsilon}{2}+\delta\right)^2+\frac{\varepsilon\pi^2}{2}\left(1+\frac{\varepsilon}{2}\right)+\frac{3\left(\frac{\varepsilon}{2}+\delta\right)}{1-\left(\frac{\varepsilon}{2}+\delta\right)} \leq \frac{1}{2}
\end{equation}
\begin{eqnarray}\label{epsdel-ineq3}
\frac{1}{24\pi^2} & \geq & R_1(\varepsilon,\delta):=\varepsilon\left(B_{221}+\frac{1}{8}B_{200}\right)+\frac{1}{16}\left(\frac{\varepsilon}{2}+\delta\right)^2\\
&+& 2\pi B_{111}\left( \left(\frac{\varepsilon}{2}+\delta \right)\left( 1+\frac{\varepsilon}{2}\right)+\varepsilon B_{210} \right)\nonumber\\
&+& \frac{\frac{\varepsilon}{2}+\delta}{1-\left(\frac{\varepsilon}{2}+\delta\right)}\left[ B_{120}+\frac{2}{\pi^2\left( 1-\frac{\varepsilon}{2}-\delta\right)^2}\left( 1+\varepsilon\pi^2B_{200}+\frac{\frac{\varepsilon}{2}+\delta}{1-\left(\frac{\varepsilon}{2}+\delta\right)}\right)\right]\nonumber
\end{eqnarray}
\begin{equation}\label{epsdel-ineq3bis}
\frac{1}{16\sqrt{3}(\pi^2+5)}\geq\sqrt{\frac{\frac{\varepsilon}{2}+\delta}{1-\left(\frac{\varepsilon}{2}+\delta\right)}}
\end{equation}
\begin{equation}\label{epsdel-ineq4}
\frac{1}{24\pi^2}  \geq  R_2(\varepsilon,\delta):=R_1(\varepsilon,\delta)+8B_{111}\left( \left(\frac{\varepsilon}{2}+\delta\right)\left(1+\frac{\varepsilon}{2}\right)\pi + \varepsilon\pi B_{210}\right)
\end{equation}
\begin{equation}\label{epsdel-ineq4bis}
\varepsilon\leq\frac{1}{\pi^2 B_{211}}\left(\frac{1}{\pi\sqrt{192}}-\frac{1}{3852}\right) 
\end{equation}
\begin{equation}\label{epsdel-ineq5}
\frac{1}{15408\pi}\geq \varepsilon\pi B_{211}+\frac{\frac{\varepsilon}{2}}{\pi\left(1-\left(\frac{\varepsilon}{2}+\delta\right)\right)}+\left(\frac{\varepsilon}{2}+\delta\right)^2\frac{\pi}{2}
\end{equation}
\begin{equation}\label{epsdel-ineq6}
\frac{1}{2\pi^2 15408^2}\geq \frac{\frac{\varepsilon}{2}+\delta}{1-\left(\frac{\varepsilon}{2}+\delta \right)}\left(B_{120}+\frac{3}{2}B_{121}+2B_{110}B_{111}+\frac{36}{5\pi^2}\right)
\end{equation}
\begin{equation}\label{epsdel-ineq7}
\frac{1}{100\pi^2 15408^2} \geq \frac{\frac{\varepsilon}{2}+\delta}{1-\left(\frac{\varepsilon}{2}+\delta \right)}\left(B_{120}+\frac{3}{2}B_{121}+2B_{110}B_{111}+\frac{36}{5\pi^2}\right)
\end{equation}
\begin{eqnarray}\label{epsdel-ineq8}
\frac{\sqrt{3}}{8\pi^2} &\geq& \frac{\left(\frac{\varepsilon}{2}+\delta\right)}{\left(1-\frac{\varepsilon}{2}-\delta\right)^3\pi^2}\left(\frac{6}{5}+\left(1+\frac{\varepsilon}{2}\right)^2\pi^2\right)+\varepsilon\left( B_{221}+\frac{16}{5}B_{211}\right)\\
&+& 7\pi B_{111}\left(\left(\frac{\varepsilon}{2}+\delta\right)\left(1+\frac{\varepsilon}{2}\right)+\varepsilon B_{210}\right)+\frac{\left(\frac{\varepsilon}{2}+\delta\right)}{1-\left(\frac{\varepsilon}{2}+\delta\right)}\left( B_{120}+\frac{36}{5\pi^2}\right)\nonumber
\end{eqnarray}

\subsection{List of constants and conditions on $\varepsilon$ and $\delta=\delta_2$ for Section \ref{sec-zero}}\label{app-sec-zero}
\begin{equation}\label{epsdel-ineq9}
\left(1+\frac{\varepsilon}{2}\right)\delta\leq\frac{\pi}{2}
\end{equation}
$$c_{11}=\sup_{\tau\in\left[0,\frac{\pi}{2}\right]}\left| \frac{\tau^2+\tau\cos\tau\sin\tau-2\sin^2\tau}{\tau^5\sin^2\tau}-\frac{2}{45\tau}\right|$$
$$c_{12}=\sup_{\tau\in\left[0,\frac{\pi}{2}\right]}\left| \frac{2(\sin\tau-\tau\cos\tau)}{\tau^3\sin^3\tau}-\frac{2}{3\tau^3}\left(1+\frac{2\tau^2}{5}\right)\right|$$
$$c_{13}=\sup_{\tau\in\left[0,\frac{\pi}{2}\right]}\left| \frac{2(\sin\tau-\tau\cos\tau)}{\tau^5\sin\tau}-\frac{2}{3\tau^3}\left(1+\frac{\tau^2}{15}\right)\right|$$
$$c_{14}=\sup_{\tau\in\left[0,\frac{\pi}{2}\right]}\left| \frac{4(\sin^2\tau-\tau^2)}{\tau^5\sin^2\tau}+\frac{4}{3\tau^3}\left(1+\frac{\tau^2}{5}\right)\right|$$
$$c_{15}=\sup_{\tau\in\left[0,\frac{\pi}{2}\right]}\left|\frac{\tau-\sin\tau}{\tau\sin\tau}\right|,\ c_{16}=\sup_{\tau\in\left[0,\frac{\pi}{2}\right]}\left|\frac{\tau^2\cos\tau-\sin^2\tau}{\tau\sin^2\tau}\right|,$$
$$c_{17}=\sup_{\tau\in\left[0,\frac{\pi}{2}\right]}\left| \frac{2\cos\tau-2+\tau\sin\tau}{\tau^6}+\frac{1}{12\tau^2}\right|$$
\begin{equation}\label{epsdel-ineq10}
\left(1+\frac{\varepsilon}{2}\right)\delta\leq\frac{2}{\sqrt{5}}
\end{equation}
$$C_{2}=\frac{338C_1^3\pi^{11}}{8}+\frac{206}{100}\left[ 8(c_{15}c_{6}+c_{16}c_{7}) +(c_{6}+c_{7})c_{15}\right]$$
$$C_{3}=\frac{268C_1^3\pi^{11}}{8}+\frac{411}{100}\left[ 2(c_{15}c_{6}+c_{16}c_{7}) + (c_{6}+c_{7})c_{15}\right]$$
\begin{equation}\label{epsdel-ineq11}
\frac{1}{180}\geq \varepsilon\ ( C_{2}+ 19c_{17}\ \delta) + \frac{115}{100}\ \delta\ (c_{11} +c_{12} +\frac{1}{2} c_{14})
\end{equation}
\begin{equation}\label{epsdel-ineq12}
 \frac{1}{180}\geq \varepsilon\ ( C_{3}+ 13c_{17}\ \delta) + \frac{115}{100}\delta\ (c_{13} +\frac{1}{2} c_{14})
\end{equation}

\subsection{A constant and conditions on $\varepsilon$ for Section \ref{sec-else}}\label{app-sec-else}
$$
C_{4}=338\pi^{10}C_1^3+20\pi(c_6+c_7)
$$
\begin{equation}\label{eps-ineq-else1}
\varepsilon \leq 2.19\times 10^{-3}\ \frac{\sin^31}{C_4}\delta_2^2
\end{equation}
\begin{equation}\label{eps-ineq-else2}
\left(1+\frac{\varepsilon}{2} \right)\left( 1+\frac{2\varepsilon}{\delta_1\left(1-\frac{\delta_1^2\pi^2}{96} \right)}\right)\leq\frac{6}{5}
\end{equation}
\begin{equation}\label{eps-ineq-else3}
\varepsilon \leq \frac{1.3\times 10^{-4}}{\pi^3C_{4}}\ \sin^3\left(\frac{\pi}{4}\delta_1\right)
\end{equation}
\bigskip

\noindent\textbf{Acknowledgment:} We thank Alessio Figalli and Ludovic Rifford for their keen interest in a preliminary draft
of this paper (44 pages, communicated to Figalli by Ge on January 14th, 2009). The second author wishes to thank Young--Heon Kim and Robert McCann for
stimulating conversations while visiting IPAM at UCLA for a spring workshop on optimal transport (April 2008). Finally, we would like to thank the Referee for a valuable comment.
{\begin{flushright}{Philippe DELANOE\\Universit\'e de Nice--Sophia
Antipolis, Facult\'e des Sciences\\Laboratoire J.--A. Dieudonn\'e, Parc Valrose\\F-06108 Nice Cedex 2\\e-mail: {\sf
Philippe.DELANOE@unice.fr}\\
\medskip
Yuxin GE\\
Universit\'e Paris--Est Cr\'eteil Val de Marne\\
Facult\'e des Sciences et Technologie\\
Centre de Math\'ematiques\\
61 avenue du G\'en\'eral De Gaulle\\
F-94010 Cr\'eteil Cedex\\
e-mail: {\sf ge@univ-paris12.fr}}
\end{flushright}}
\end{document}